\documentclass[10pt]{amsart}
\usepackage{amsmath,amssymb,amsthm}
\newtheorem{thm}{Theorem}[section]
\newtheorem{lem}[thm]{Lemma}
\newtheorem{prop}[thm]{Proposition}

\numberwithin{equation}{section}

\theoremstyle{example}
\newtheorem{example}[thm]{Example}

\theoremstyle{definition}

\theoremstyle{notation}
\newtheorem{notation}[thm]{Notation}
\theoremstyle{remark}
\newtheorem{remark}[thm]{Remark}

\newcommand{\C}{{\mathbb C}}
\newcommand{\N}{{\mathbb N}}
\newcommand{\Z}{{\mathbb Z}}
\newcommand{\R}{{\mathbb R}}

\title
[Invariant measures for BO ]
{Gaussian measures associated to the higher order conservation laws of the  Benjamin-Ono equation }
\author[Nikolay Tzvetkov]{Nikolay~Tzvetkov}
\author[Nicola Visciglia]{Nicola Visciglia}
\address{D\'epartement de Math\'ematiques, Universit\'e de Cergy-Pontoise, 2, 
avenue Adolphe Chauvin, 95302 Cergy-Pontoise  
Cedex, France and Institut Universitaire de France}\email{nikolay.tzvetkov@u-cergy.fr}
\address{Universit\`a Degli Studi di Pisa Dipartimento di Matematica "L. Tonelli"
Largo Bruno Pontecorvo 5 I - 56127 Pisa. Italy}\email{ viscigli@dm.unipi.it}
\begin{document}
\maketitle
\begin{abstract}
Inspired by the work of Zhidkov on the KdV equation, we perform a construction of  weighted gaussian measures 
associated to the higher order conservation laws of the Benjamin-Ono equation. 
The resulting measures are supported by Sobolev spaces of increasing regularity. We also prove 
a property on the support of these measures
leading to the conjecture that they are indeed invariant by the flow of the Benjamin-Ono equation. 
\end{abstract}
%

\section{Introduction and statement of the results}
\subsection{Measures construction} 
The main goal of this article is to construct weighted gaussian measures associated with an arbitrary conservation law of the 
Benjamin-Ono equation (BO), and thus to extend the result of the first author \cite{tz} which deals only with the first conservation law.
The analysis contains several significant elaborations with respect to \cite{tz}, it requires an understanding of the interplay between the structure of the conservation laws of
the Benjamin-Ono equation and the probabilistic arguments involved in the renormalization procedure defining the measures. 

Let us recall that just like the KdV equation, the Benjamin-Ono equation is a basic dispersive PDE describing the propagation of one directional, long, small amplitude waves. The difference between the KdV and BO equations is that the KdV equation describes surface waves while the Benjamin-Ono equation models the propagation of internal waves.
These models have rich mathematical structure from both the algebraic and analytical viewpoints. 
In particular they have an infinite sequence of conservation laws. These aspects will be heavily exploited in the present work.  
 
Consider now the Benjamin-Ono equation
\begin{equation}\label{bo}
\partial_t u + H\partial_x^2 u + u\partial_x u=0,
\end{equation}
with periodic boundary conditions (for simplicity throughout the paper we fix the period to be equal to $2\pi$).
In \eqref{bo}, $H$ denotes the Hilbert transform acting on periodic distributions.
Thanks to the work of Molinet \cite{M} \eqref{bo} is globally well-posed in $H^s$, $s\geq 0$ 
(see \cite{T, IK, BP} for related results in the case when \eqref{bo} is posed on the real line).

It is well-known that (smooth) solutions to \eqref{bo} satisfy infinite number of 
conservation laws (see e.g. \cite{Mats, ABFS}).
More precisely for $k\geq 0$ an integer, there is a conservation law of \eqref{bo} of the form 
\begin{equation}\label{strucRmezzi}
E_{k/2}(u)= \|u\|_{\dot{H}^{k/2}}^2 + R_{k/2} (u) 
\end{equation}
where $\dot{H}^s$ denotes the homogeneous Sobolev norm on periodic functions, and all the terms that appear in $R_{k/2}$ 
are homogeneous of the order larger or equal than three in $u$.
In Section~\ref{conserlawsbo}, we will describe in more details the structure of $R_{k/2}$
for large $k$. Next we explicitly write the conservation laws
$E_{k/2}$ for $k=0,1,2,3,4$:
\begin{eqnarray*}
E_0(u) & = & \|u\|_{L^2}^2;
\\
E_{1/2}(u) & = & \|u\|_{\dot{H}^{1/2}}^2 + \frac 13 \int u^3 dx;
\\
E_1(u) & = & \|u\|_{\dot{H}^1}^2 + \frac 34  \int u^2 H(u_x) dx + \frac 18 \int u^4 dx;
\\
E_{3/2}(u) & = & \|u\|_{\dot{H}^{3/2}}^2 - \int [\frac 32 u (u_x)^2 + \frac 12 u H(u_x)^2] dx 
\\
& &
- \int [\frac 13 u^3 H (u_x) + \frac 14 u^2 H(u u_x)] dx-\frac 1{20} \int u^5dx;
\\
E_2(u) & = & \|u\|_{\dot{H}^2}^2 -\frac 54 \int  [(u_x)^2 H u_x  + 2 u u_{xx} Hu_x] dx
\\
& &
 + \frac 5{16}\int  [5 u^2 (u_x)^2 + u^2 H(u_x)^2+ 2u H(\partial_x u) H (u u_x) ] dx
\\
& &
+ \int [\frac 5{32} u^4 H(u_x) + \frac 5{24} u^3 H(u u_x)] 
dx+ \frac 1{48} \int u^6 dx
\end{eqnarray*}
where $\int$ is understood as the integral on the period $(0, 2\pi)$.

Following the work by Zhidkov \cite{zh} (see also \cite{B, LRS}),
one may try to define an invariant measure for \eqref{bo} 
by re-normalizing the formal measure $\exp(-E_{k/2}(u))du$.
This re-normalization is a delicate procedure. One possibility would be first to re-normalize 
$\exp(-\|u\|_{\dot{H}^{k/2}}^2)du$ as a gaussian measure on an infinite dimensional space and then to show that the factor 
 $\exp(-R_{k/2}(u))$ is integrable with respect to this measure.
 
Since  $\exp(-\|u\|_{\dot{H}^{k/2}}^2)$ factorizes as an infinite product when we express $u$ as a Fourier series, we can 
define the re-normalization of $\exp(-\|u\|_{\dot{H}^{k/2}}^2)du$ as the gaussian measure induced 
by the random Fourier series 
\begin{equation}\label{randomized}
\varphi_{k/2}(x, \omega)=\sum_{n\neq 0} 
\frac{\varphi_n(\omega)}{|n|^{k/2}} e^{{\bf i}nx}
\end{equation}
(one may ignore the zero Fourier mode since the mean of $u$ is conserved by the flow of \eqref{bo}).
In \eqref{randomized}, $(\varphi_{n}(\omega))_{n\neq 0}$ is a sequence of standard complex gaussian variables 
defined on a probability space $(\Omega, {\mathcal A}, p)$ such that $\varphi_{n}=\overline{\varphi_{-n}}$
(since the solutions of \eqref{bo} should be real valued) and  
$(\varphi_{n}(\omega))_{n>0}$ are independent.
Let us denote by $\mu_{k/2}$ the measure induced by \eqref{randomized}. 
One may easily check that $\mu_{k/2}(H^s)=1$ 
for every $s<(k-1)/2$ while
$\mu_{k/2}(H^{(k-1)/2})=0$. 

In view of the previous discussion, one may consider $\exp(-R_{k/2}(u))d\mu_{k/2}$ as a candidate 
of invariant measure for \eqref{bo}.
There are two obstructions to do that, the first one already 
appears in previous works on the NLS equation (see \cite{B, LRS}) 
and the KdV equation (see \cite{zh}), while 
the second one is specific to the Benjamin-Ono equation.
The first obstruction is that $ \exp(-R_{k/2}(u))$ is not integrable with respect to $d\mu_{k/2}(u)$. 
This problem may be 
resolved by restricting to invariant sets, which means to replace $ \exp(-R_{k/2}(u))$  by 
\begin{equation}\label{restrict}
\prod_{j=0}^{k-1} 
\chi_R (E_{j/2}(u)) 
e^{-R_{k/2}(u)}\,\,,
\end{equation}
where $\chi_{R}$ is a cut-off function defined as $\chi_R(x)=\chi (x/R)$
with $\chi:\R \rightarrow \R$ a continuous,
compactly supported function such that $\chi(x)=1$ for every $|x|<1$. 
In the context of KdV or NLS, the function defined in \eqref{restrict} is 
integrable with respect to the corresponding gaussian measure.
Moreover if one takes the reunion over $R>0$ of the supports of the functions 
\eqref{restrict}, then one obtains a set containing the support of
$\mu_{k/2}$. However, in the context of the Benjamin-Ono 
equation, the restriction to invariant sets does not work as in \eqref{restrict} 
because for every $R$ the following occurs: $\chi_{R}(E_{(k-1)/2}(u))=0$ almost surely on the support of $\mu_{k/2}$.
One of the main points of this paper is to resolve this difficulty.
This will be possible since one controls the way that $E_{(k-1)/2}(u)$ diverges on the support of $\mu_{k/2}$.
More precisely, for $N\geq 1$ and $k\geq 2$, we introduce the function
\begin{equation}\label{Femme+1}
F_{k/2,N, R}(u)
=\Big(\prod_{j=0}^{k-2} \chi_R (E_{j/2}(\pi_N u)) \Big)
\chi_R (E_{(k-1)/2}(\pi_N u)-\alpha_N) 
e^{-R_{k/2}(\pi_N u)}
\end{equation}
where
$\alpha_N=\sum_{n=1}^N \frac 1n$ and $\pi_N$ is the Dirichlet projector 
on Fourier modes $n$ such that $|n|\leq N$.
Here is our first result.
\begin{thm}\label{main}
For every $k\in \N$ with $k\geq 2$, there exists a $\mu_{k/2}$ measurable function 
$F_{k/2,R }(u)$ such that $F_{k/2,N, R}(u)$ converges to $F_{k/2,R }(u)$ in $L^q(d\mu_{k/2})$ 
for every $1\leq q<\infty$.
In particular $F_{k/2,R }(u)\in L^q(d\mu_{k/2})$. 
Moreover, if we set $d\rho_{k/2,R}\equiv F_{k/2,R }(u)d\mu_{k/2}$, we have
$$
\bigcup_{R>0}{\rm supp}(\rho_{k/2,R})={\rm supp}(\mu_{k/2}).
$$
\end{thm}
The above result for $k=1$ was obtained by the first author in \cite{tz}.
Many of the probabilistic techniques involved in the proof of Theorem \ref{main}
are inspired by \cite{B2}. We also refer to \cite{BS} where 
in the context of the $2d$ NLS the authors use the Wick ordered
$L^2$-cutoff, i.e. a truncation of the $L^2$-norm that depends on
the parameter $N$.

We conjecture that the measures $\rho_{k/2, R}$, $k=2,3,\cdots$ constructed in Theorem~\ref{main}
are invariant by the flow of the Benjamin-Ono equation established by Molinet \cite{M}, 
at least for even values of $k$.
In the sequel, for shortness, we denote $\rho_{k/2, R}$ by $\rho_{k/2}$.
\subsection{A property on the support of the measures}
Let us now give our argument in support of the above-stated conjecture. 
For $N\geq 1$, we introduce the truncated Benjamin-Ono equation:
\begin{equation}\label{truncBO}
\partial_t u + H\partial_x^2 u + \pi_N \big((\pi_{N}u)\partial_x(\pi_{N} u)\big)=0.
\end{equation}
As in \cite{BTT}, one can define a global solution of \eqref{truncBO} 
for every initial data $u(0)\in L^2(S^1)$.
Indeed, one obtains that $(1-\pi_{N})u(t)$ is given by the free Benjamin-Ono 
evolution with data $(1-\pi_{N})u(0)$, while
$\pi_{N}u(t)$ evolves under an $N$-dimensional ODE. This ODE has a well-defined global dynamics since 
the $L^2$ norm is preserved. 

The main problem that appears when one tries 
to prove the invariance of $\rho_{k/2}$ is that even if $E_{k/2}$ 
are invariants for the Benjamin-Ono equation they are 
not invariant under \eqref{truncBO}.
The invariance, however, holds in a suitable asymptotic sense as we explain below. 
Let us introduce the  real-valued function
$G_{k/2,N}$, measuring the lack of conservation of $E_{k/2}$ under the truncated flow \eqref{truncBO}, via the following relation
\begin{equation}\label{G}
\frac{d}{dt} E_{k/2}(\pi_{N}u(t))=G_{k/2,N}(\pi_{N}u(t)),
\end{equation}
where $u(t)$ solves \eqref{truncBO}.
\\
Denote by 
$\Phi_{N}$ the flow of \eqref{truncBO} and set $d\rho_{N}(u)\equiv F_{k/2,N,R}(u)d\mu_{k/2}(u)$ 
so that by Theorem~\ref{main}, $\rho_{N}$ converges in a strong sense to $\rho_{k/2}$
(the densities converge in any $L^p(d\mu_{k/2})$, $p<\infty$).
By using the Liouville theorem, one shows that for every $\mu_{k/2}$ measurable set $A$,
$$
\rho_{N}(\Phi_{N}(t)(A))=
\int_{A}
e^{-\int_{0}^{t}G_{k/2,N}(\pi_{N}\Phi_N(\tau)(u(0))d\tau}d\rho_{N}(u(0))
+o(1)\,.
$$
Hence, a main step towards a proof of the invariance of $\rho_{k/2}$ is to show that
\begin{equation}\label{object}
\int_{0}^{t}G_{k/2,N}(\pi_{N}u(\tau))d\tau
\end{equation}
converges to zero,  where $u(\tau)$ is a solution of \eqref{truncBO}, 
with $u(0)$ on the support of  $\mu_{k/2}$. Such a property is relatively easy to be established if $u(0)$ has slightly more regularity
than the typical Sobolev regularity on the support of $\mu_{k/2}$.
At the present moment, we are not able to prove such a property on the support of $\mu_{k/2}$.
We shall, however, prove it if we make a first approximation which consist of 
replacing $u(\tau)$ by $u(0)$ in \eqref{object}.
Here is the precise statement.
\begin{thm}
\label{invariance}
For every $k\geq 6$ an even integer, we have
$$
\lim_{N\rightarrow \infty}\|G_{k/2,N}(\pi_N u)\|_{L^q(d\mu_{k/2})}=0,\hbox{ } \forall q\in [1, \infty),
$$
where $G_{k/2,N}$ is defined by \eqref{G}.
\end{thm}
Let us remark that the lack of invariance of conservation laws for the corresponding truncated flows is a problem that appears
also in other contexts. We refer in particular to the papers \cite{NORS} and \cite{zh}, where this difficulty is resolved in the cases of the DNLS and KdV equations respectively.
\subsection{Comparison with the KdV equation}\label{kdvrie}
Next we explain why the measures construction in the context of the Benjamin-Ono equation is much more involved compared with the case of the KdV equation. The main difference is that
in the KdV equation the dispersion (of the linear part) is of lower order compared with the Benjamin-Ono equation. This fact makes the perturbative treatment of the 
nonlinearity more complicated. Let us recall that a similar observation applies to the Cauchy problem analysis (see \cite{BP, IK, M, T}).

Now, we recall the approach of Zhidkov (see \cite{zh})
to prove the existence of invariant measures associated with the periodic KdV equation
\begin{equation}\label{KdV}
\partial_t u + \partial_x^3 u + u \partial_x u=0.\end{equation}
This equation has a rich
structure from both the algebraic and analytic
viewpoint. In particular the solutions to \eqref{KdV} have an infinite sequence of conversation laws.
More precisely, for every $m\geq 0$ there exists a polynomial  
$$p_m(v, \partial_x v,..., \partial_x^m v )$$ such that
$\frac d{dt} \int p_m (u(t,x), \partial_x u(t, x),..., \partial_x^m u(t, x))dx=0$,
provided that $u$ is a solution of \eqref{KdV},
where $\int...dx$ denotes the integral on the period.
More precisely,  the conservation laws have  the following structure 
$$\|u\|_{\dot{H}^m}^2 +  \int q_m( u,.., \partial_x^{m-1} u) dx  .$$ 
By using the Sobolev embedding $H^1\subset L^\infty$, it is easy to check that the function
$$\Big(\prod_{j=0}^{m-1} \chi_R (E_j(u)) \Big) e^{-\int 
q_m(u,.., \partial_x^{m-1} u) dx },$$ 
where $\chi_R$ is defined as in \eqref{restrict},
is not trivial and belongs to the space
$L^\infty(d\mu_m)$ (where $\mu_m$ is the Gaussian measure induced by 
\eqref{randomized} for $k=2m$), provided that $m$ is large enough. In particular
the measure
\begin{equation}\label{candKdV}
\Big(\prod_{j=0}^{m-1} \chi_R (E_j(u)) \Big) e^{-\int 
q_m(x, u,.., \partial_x^{m-1} v)dx } d\mu_m
\end{equation} 
is a meaningful non-trivial candidate for an invariant measure.\\
In order to prove the invariance of the above measure 
we introduce, following \cite{zh},
a family of truncated problems
\begin{equation}\label{truncKdV}
\partial_t u + \partial_x^3 u + \pi_N \big((\pi_{N}u)\partial_x(\pi_{N} u)\big)=0,
\end{equation}
where $\pi_N$ is the Dirichlet projector on the $n$ Fourier modes such that $|n|\leq N$.
Once again, the main difficulty is because of the fact that the quantity
$$\int p_m(u(t,x), \partial_x u(t,x),..., \partial_x^m u(t, x)) dx$$
is no longer invariant along the flow of the truncated problem \eqref{truncKdV}.
However if 
$u_0\in H^{m-1}$ 
then 
\begin{equation}\label{asymptlw}
\lim_{N\rightarrow \infty} \frac d{dt} \int p_m(\pi_N u(t,x), \partial_x \pi_N u(t,x),..., 
\partial_x^m \pi_N u(t, x)) dx=0,
\end{equation}
where $u(t,x)$ are solutions
to \eqref{truncKdV} with initial data $u_{0}$.
Roughly speaking \eqref{asymptlw} means that the quantities $$\int p_m (\pi_N u(t,x), \partial_x 
\pi_N u(t,x),..., \partial_x^m \pi_N u(t, x)) dx$$
are asymptotically in $N$ almost conservation laws
for solutions to the truncated flow. In particular for large $N$,  
the classical finite dimensional Liouville invariance theorem  turns out 
to be almost true for the flow associated with \eqref{truncKdV},
and it allows us to conclude the proof of the invariance
of \eqref{candKdV} along the flow associated with
\eqref{KdV} via a limit argument.\\
Hence the main point is to prove \eqref{asymptlw}.
Following Zhidkov (see Lemma IV.3.5, page 127 in \cite{zh}) 
there is an explicit formula to compute
the expression on the l.h.s. in \eqref{asymptlw}.
More precisely if $u$ solves \eqref{truncKdV}
then
\begin{equation}\label{zhid}
\frac d{dt} \int p_m(\pi_N u(t,x), \partial_x \pi_N u(t,x),..., 
\partial_x^m \pi_N u(t, x)) dx\end{equation}
$$=\sum_{j=0}^m
\int \Big( \frac{\partial p_m }{\partial_{x}^j u}\Big)_{|\partial_{x}^ju= 
\partial_{x}^j \pi_{>N} ((\pi_N u) \partial_x (\pi_N u)), \partial_x^k u= 
\partial_x^k \pi_N u \hbox{ for } k\neq j} dx,
$$
where $\pi_{>N}$ is the projector on the $n$ Fourier modes such that $|n|>N$
(for an explanation of formula \eqref{zhid} see Section \ref{pNstar}). 
It is easy to see that the most delicate term that appears in the r.h.s. above is the one coming
from the cubic part of the conservation laws, i.e.
$\int u (\partial_x^{m-1} u)^2 dx$.
More precisely we have to estimate the following term
\begin{equation}\label{kdvworst}
\int (\pi_N u) \partial_x^{m-1}(\pi_N u) \partial_x^{m-1} \pi_{>N} ((\pi_N u) \partial_x 
(\pi_N u)) dx,
\end{equation}
as $N\rightarrow \infty$.
Notice that after developing
the $(m-1)$-derivative of the product, we get 
an integral that involves the product of $\partial_x^{m-1} (\pi_N u)$ 
and $\partial_x^{m}(\pi_N u)$,
and hence after a fractional integration by parts we get 
a derivative of order $m-1/2$. This is the main source of difficulty since the Gaussian measure
$d \mu_m$ is supported on the Sobolev spaces $H^{m-1/2 -\epsilon}$
for any $\epsilon>0$. 
This problem is solved by Zhidkov by using
a clever integration by parts. 
Indeed, if we develop the
$(m-1)$-derivative of the product in \eqref{kdvworst}, using the Leibnitz rule, 
we get the following (bad) term
$$\int (\pi_N u) \partial_x^{m-1}(\pi_N u) \pi_{>N} 
((\pi_N u) \partial_x^m (\pi_N u)) dx$$
$$= \int \pi_{>N} ((\pi_N u) \partial_x^{m-1}(\pi_N u))  
\partial_x (\pi_{>N}((\pi_N u) \partial_x^{m-1} (\pi_N u))) dx$$
$$-\int \pi_{>N} ((\pi_N u) \partial_x^{m-1}(\pi_N u)) \pi_{>N} 
((\pi_N \partial_x u) \partial_x^{m-1} (\pi_N u)) dx.$$
The worst term in the r.h.s. seems to be the first one,
since it involves the product of a derivative of order $m$ and a derivative of order $m-1$. However
this term is zero since
it can be written as follows:
$$\frac 12 \int \partial_x (\pi_{>N}(\pi_N u) \partial_x^{m-1} (\pi_N u))^2dx=0.$$

By looking at the structure of the conservation laws of the Benjamin-Ono
equation, it is easy to check that the situation is a priori much worse.
In fact if $E_m$ is the conservation law (for the Benjamin-Ono equation) with leading term
$\|u\|_{H^m}^2$, and if we repeat the same construction 
as in \eqref{zhid} (where $p_m$ is replaced by the density of $E_m$),
then the cubic part 
produces a contribution that involves a derivative of order $m$
which is very delicate since the Gaussian measure 
$d\mu_m$ is supported on $H^{m-1/2 - \epsilon}$. Moreover, in $E_m$ 
the terms homogeneous of order four  
involve a derivative of order $m-1/2$ (this difficulty can be compared with the one we 
met above to treat the contribution coming from the cubic part of the conservation laws of
KdV).
The second main result of this paper (Theorem~\ref{invariance}) is essentially 
saying that we are able to find a key cancellation which eliminates the terms containing (at first glance) too many derivatives.
We believe that this result is of independent interest and that it will play a role in the future analysis on the issues considered here.
\\

Next we fix some notations.

\begin{notation}
We shall denote by $H^s$ (and in some cases $H^s_x$)
the Sobolev
spaces of $2\pi$-periodic functions;\\
$L^p$ (and in some cases $L^p_x$) is 
the $2\pi$-periodic Lebesgue space;\\
$L^q_\omega$ is the Lebesgue 
space with respect to the probability measure 
$(\Omega, {\mathcal A}, p)$, which in turn is the domain of definition
of the random variables 
$\varphi_n(\omega)$ in \eqref{randomized};\\
if $f(x)$ is  a $2\pi$-periodic function then
$\int f(x) dx= \int_0^{2\pi} f(x) dx$;\\
the operator $H$ is the usual Hilbert transform
acting on $2\pi$-periodic functions;\\
for every $N\in \N$ the constant $\alpha_N$ is equal to
$\sum_{n=1}^N \frac 1n$;\\
for every $k\in\ N$ the quantities
$E_{k/2}$ and $R_{k/2}$
are related as in \eqref{strucRmezzi}, where $E_{k/2}$
is a conservation law for the Benjamin-Ono equation.  
\\
Some other notations will be fixed in Section~\ref{conserlawsbo}.
\end{notation}
The remaining part of the paper is devoted to the proof of Theorems~\ref{main},  \ref{invariance}.
\section{On the structure of the conservation laws of the Benjamin-Ono equation}\label{conserlawsbo}
In this section, we describe the form of the Benjamin-Ono equation conservation laws which is suitable for the proof of the our results announced in the introduction.
Our reference in this discussion is the book by Matsuno \cite{Mats}. 

We now fix some notations.
Given any function $u(x)\in C^\infty(S^1)$, we define
\begin{eqnarray*}
{\mathcal P}_1(u) & = & \{\partial_x^{\alpha_1} u, 
H\partial_x^{\alpha_1} u|\alpha_1\in \N\},
\\
{\mathcal P}_2(u) & = & \{\partial_x^{\alpha_1} 
u\partial_x^{\alpha_2} u, 
(H \partial_x^{\alpha_1} u)\partial_x^{\alpha_2} u,
(H \partial_x^{\alpha_1} u)(H \partial_x^{\alpha_2} u)|\alpha_1,\alpha_2\in \N\}
\end{eqnarray*}
and in general by induction
\begin{multline*}
{\mathcal P}_n(u)=\Big \{\prod_{l=1}^k H^{i_l}p_{j_l}(u)|
i_1,...,i_k\in \{0,1\}, 
\\
\sum_{l=1}^k j_l=n, k\in \{2,...,n\}
\hbox{ and } p_{j_l}(u)\in {\mathcal P}_{j_l}(u)\Big \},
\end{multline*}
where $H$ is again the Hilbert transform.
\begin{example}
The elements belonging to 
${\mathcal P}_3(u)$ are the following ones: 
$$\partial_x^{\alpha_1}u
\partial_x^{\alpha_2}u\partial_x^{\alpha_3}u,
\partial_x^{\alpha_1}u\partial_x^{\alpha_2}u
(H \partial_x^{\alpha_3}u),
\partial_x^{\alpha_1}u H 
(\partial_x^{\alpha_2}u\partial_x^{\alpha_3}u),$$$$ 
\partial_x^{\alpha_1} u(H \partial_x^{\alpha_2} u)
(H \partial_x^{\alpha_3}u),
\partial_x^{\alpha_1}uH (\partial_x^{\alpha_2}
u(H \partial_x^{\alpha_3}u)),
H \partial_x^{\alpha_1}u H (\partial_x^{\alpha_2}
u\partial_x^{\alpha_3} u),$$
$$(H \partial_x^{\alpha_1}u)(H \partial_x^{\alpha_2}u)
(H \partial_x^{\alpha_3}u),
(H \partial_x^{\alpha_1}u) H (\partial_x^{\alpha_2}u
(H\partial_x^{\alpha_3}u)),
\partial_x^{\alpha_1}u H((H \partial_x^{\alpha_2}u)
(H\partial_x^{\alpha_3}u))$$$$\hbox{ where } \alpha_1, \alpha_2, \alpha_3\in \N.$$
\end{example}
\begin{remark}
Roughly speaking an element in 
${\mathcal P}_n(u)$ involves the product of 
$n$ derivatives $\partial_x^{\alpha_1} u, .., 
\partial^{\alpha_n}_x u$ 
in combination with the Hilbert transform $H$ 
(that can appear essentially in an arbitrary way in front of the factors and eventually in front of a group of factors).
\end{remark}
Notice that for every $n$ the simplest element belonging 
to $\mathcal P_n(u)$ has the following structure:
\begin{equation}\label{ptilde}
\prod_{i=1}^n \partial_x^{\alpha_i} u, \alpha_i\in \N.
\end{equation}
In particular we can define the map
$${\mathcal P}_n(u)\ni p_n(u)\rightarrow 
\tilde p_n(u)\in {\mathcal P_n}(u)$$
that associates to every $p_n(u) \in {\mathcal P}_n(u)$ the unique element
$\tilde p_n(u)\in {\mathcal P}_n(u)$ having the structure given in \eqref{ptilde}
where $\partial_x^{\alpha_1}u, \partial_x^{\alpha_2}u,
..., \partial_x^{\alpha_n}u$
are the derivatives involved in the expression of $p_n(u)$
(equivalently $\tilde p_n(u)$ is obtained from $p_n(u)$ by 
erasing all the Hilbert transforms $H$ that appear in $p_n(u)$).\\
Next, we associate to every $p_n(u)\in {\mathcal P_n}(u)$ 
two integers as follows:\\
$$\hbox{ if }\tilde p_n(u)=\prod_{i=1}^n \partial_x^{\alpha_i} u
\hbox{ then }$$ 
\begin{equation}\label{norm}|p_n(u)|:=\sup_{i=1,..,n} \alpha_i
\end{equation}
and
\begin{equation}\label{doublenorm}\|p_n(u)\|:=\sum_{i=1}^n \alpha_i.\end{equation}
We are ready to describe the structure of the 
conservation laws satisfied by the
Benjamin-Ono equation.
Given any even $k\in \N$, i.e. $k=2n$, the energy $E_{k/2}$ 
has the following structure:

\begin{equation}\label{even}E_{k/2}(u)
=  \|u\|_{\dot{H}^{k/2}}^2 + 
\sum_{\substack{p(u)\in {\mathcal P}_3(u) s.t. \\ 
\tilde p(u)=u\partial_x^{n-1} u \partial_x^{n}u}}
c_k(p) \int p(u)dx\end{equation}
$$+
\sum_{\substack{p(u)\in {\mathcal P}_{j}(u) s.t. j=3,..., 2n+2\\ 
\|p(u)\|= 2n-j+2\\ |p(u)|\leq n-1}} c_k(p) \int p(u)dx 
$$
where $c_k(p)\in \R$ are suitable real numbers.\\
Similarly in the case of odd $k\in \N$, i.e. $k=2n+1$, 
the energy $E_{k/2}$ has the following structure:
\begin{equation}\label{odd}E_{k/2}(u)
=  \|u\|_{\dot{H}^{k/2}}^2 
+ 
\sum_{\substack{p(u)\in {\mathcal P}_3(u) s.t. \\ 
\tilde p(u)=u\partial_x^{n} u \partial_x^{n}u}}
c_k(p) \int p(u)dx\end{equation}
$$
+\sum_{\substack{p(u)\in {\mathcal P}_3(u) s.t. \\ 
\tilde p(u)=\partial_x u\partial_x^{n-1} u \partial_x^{n}u}}
c_k(p) \int p(u)dx
$$
$$+
\sum_{\substack{p(u)\in {\mathcal P}_4(u) s.t. \\ 
\tilde p(u)=u^2\partial_x^{n-1} u \partial_x^{n}u}}
c_k(p) \int p(u)dx
$$
$$+
\sum_{\substack{p(u)\in {\mathcal P}_j(u) s.t.j=3,...,2n+3\\ 
\|p(u)\|=2n-j+3\\|p(u)|\leq n-1}} c_k(p) \int p(u)dx$$
where $c_k(p)\in \R$ are suitable real numbers.\\
\begin{remark}
The expressions above should be compared with the explicit structure
of $E_{k/2}$ for $k=0,1,2,3,4$ (see the introduction).
\end{remark}
\section{Preliminary Estimates}
Along this section we shall use the notations ${\mathcal P}_j(u)$,
$p_j(u)$, $\tilde p_j(u)$ introduced in Section \ref{conserlawsbo}.
We also recall that $E_{k/2}$ denotes the conservation law
whose structure is described in \eqref{even} and \eqref{odd}
(respectively depending on the eveness or oddness of $k$).
The main result of this section is the following proposition
that will be very useful to prove Theorem \ref{main} for $k$ an even number.
\begin{prop}\label{carpi}
Let $k>0$ be a given integer. Then for every $R_1, R_2>0$ 
there is $C=C(R_1,R_2)>0$ such that
\begin{equation}\label{propgian}
\bigcap_{j=0}^{2k} \{u\in H^k| |E_{j/2}(u)|< R_1\} 
\cap \{u\in H^k| |E_{k+1/2}(\pi_N u )-\alpha_N|<R_2\} 
\end{equation}
$$\subset \{ u\in H^k | \|u\|_{H^{k}}<C\}\cap 
\{u\in H^k | \|\pi_N u\|_{\dot H^{k+1/2}}^2 - \alpha_N|<C\}, \hbox{ }
\forall N\in \N.$$
\end{prop}

\begin{remark}
The proposition above (where we choose $k=m$) implies 
that the support of the functions $F_{m+1,N, R}$ defined in 
\eqref{Femme+1} are contained in a ball of $H^m$ intersected 
with the region  $\{u\in H^m | 
\|\pi_N u\|_{\dot H^{m+1/2}}^2 - \alpha_N|<C\}$ (at least in the case $m>0$).
\end{remark}

\begin{lem}\label{LP}
For every integer $m\geq 0$ there exists $C=C(m)>0$ such that 
\begin{equation}\label{edno}
\Big|\int u\,\partial_{x}^mv\,\partial_{x}^{m+1}w \hbox{ } dx\Big|
\leq C (\|u\|_{L^{\infty}}\|v\|_{H^{m+1/2}}\|w\|_{H^{m+1/2}}
\end{equation}
$$+
\|v\|_{L^{\infty}}\|u\|_{H^{m+1/2}}\|w\|_{H^{m+1/2}}
+
\|w\|_{L^{\infty}}\|u\|_{H^{m+1/2}}\|v\|_{H^{m+1/2}})\,.
$$
\end{lem}
{\bf Proof.} We consider a Littlewood-Paley partition of unity
$1=\sum_{N}\Delta_{N}$, where $N$ takes the dyadic values, i.e. $N=2^j$, $j=0,1,2,\cdots$.
We denote by $S_N$ the operator $\sum_{N_1\leq N}\Delta_{N_1}$.
In order to prove \eqref{edno}, one needs to evaluate the expression
\begin{equation}\label{zvezda}
\sum_{N_1,N_2,N_3}
\int \Delta_{N_1}u\,
\partial_{x}^{m}\Delta_{N_2}v\, 
\partial_{x}^{m+1}\Delta_{N_3}w \hbox{ } dx\,.
\end{equation}
We consider three cases by distinguishing which is the smallest of $N_1$, $N_2$ and $N_3$.
In the sequel we shall denote by $c, C>0$ constants
that can change at each step.\\
Denote by $J_1$, the contribution of $N_1\leq \min(N_2,N_3)$ to \eqref{zvezda}.
Then $N_2\sim N_3$ and
\begin{eqnarray*}
J_1 \leq C \sum_{N_2\sim N_3}
\Big|\int
S_{c\min(N_2,N_3)}u\,
\partial_{x}^{m}\Delta_{N_2}v\, 
\partial_{x}^{m+1}\Delta_{N_3}w \hbox{ } dx 
\Big|
\\
\leq C
\sum_{N_2\sim N_3}\|u\|_{L^{\infty}}
N_2^{m}\|\Delta_{N_2}v\|_{L^2}N_{3}^{m+1}\|\Delta_{N_3} w\|_{L^2}
\\
\leq C
\|u\|_{L^{\infty}}\|v\|_{H^{m+1/2}}\|w\|_{H^{m+1/2}},
\end{eqnarray*} 
where in the last line we used the Cauchy-Schwarz inequality.
Next denote by $J_2$ the contribution of $N_2\leq \min(N_1,N_3)$ to \eqref{zvezda}.
Then
\begin{eqnarray*}
J_2 \leq C\sum_{N_1\sim N_3}
\Big|\int
\Delta_{N_1}u\,
\partial_{x}^{m}S_{c\min(N_1,N_3)}v\, 
\partial_{x}^{m+1}\Delta_{N_3}w \hbox{ }
dx \Big|
\\
\leq C
\sum_{N_1\sim N_3}\|\Delta_{N_1}u\|_{L^{2}}
(\min(N_1,N_3))^{m}\|v\|_{L^\infty}N_{3}^{m+1}\|\Delta_{N_3} w\|_{L^2}
\\
\leq C
\|v\|_{L^{\infty}}\|u\|_{H^{m+1/2}}\|w\|_{H^{m+1/2}},
\end{eqnarray*} 
Finally, we denote by $J_3$, the contribution of $N_3\leq \min(N_1,N_2)$ to \eqref{zvezda}.
Then
\begin{eqnarray*}
J_3 \leq C\sum_{N_1\sim N_2}
\Big|\int
\Delta_{N_1}u\,
\partial_{x}^{m}\Delta_{N_2}v\, 
\partial_{x}^{m+1}S_{c\min(N_1,N_2)}w
\hbox{ } dx \Big|
\\
\leq C
\sum_{N_1\sim N_2}\|\Delta_{N_1}u\|_{L^{2}}
N_2^{m}\|\Delta_{N_2}v\|_{L^2}
(\min(N_1,N_2))^{m+1}\|w\|_{L^\infty}
\\
\leq C 
\|w\|_{L^{\infty}}\|u\|_{H^{m+1/2}}\|v\|_{H^{m+1/2}},
\end{eqnarray*} 
This completes the proof of Lemma~\ref{LP}.

\hfill$\Box$

As a consequence of Lemma \ref{LP} we get the following 
useful result.
\begin{lem}\label{cubichard}
Let $m\geq 0$ be an integer and $p_3 (u)\in {\mathcal P}_3(u)$ be such that
\begin{equation}\label{strpeg}
\tilde p_3 (u)= u \partial_x^m u \partial_x^{m+1}u.
\end{equation}
Then for every $\epsilon>0$, $1<p<\infty$ such that
$\epsilon p>1$, there exists $C=C(\epsilon, p)>0$ such that:
\begin{equation}\label{hard}
\Big |\int p_3(u) dx \Big | \leq C \|u\|_{H^{m+1/2}}^2\|u\|_{W^{\epsilon,p}}.
\end{equation}
\end{lem}
{\bf Proof.}
Looking at the structure of the elements in ${\mathcal P}_3(u)$
and since we are assuming \eqref{strpeg} we can deduce by 
Lemma \ref{LP} the following estimate:
$$\Big |\int p_3(u) dx\Big |\leq C \Big( \max \{\|u\|_{H^{m+1/2}}, 
\|H u\|_{H^{m+1/2}}\}\Big)^2
\max \{\|u\|_{L^\infty}, \|Hu\|_{L^\infty}\}$$
and hence by the Sobolev embedding
$W^{\epsilon,p}\subset L^\infty$ we can continue the estimate
as follows (provided that we change the constant $C$)
$$....\leq C \|u\|_{H^{m+1/2}}^2 
\max \{\|u\|_{W^{\epsilon,p}}, \|Hu\|_{W^{\epsilon,p}}\}.$$
The proof can be completed since
the Hilbert transform $H$ is continuous
in the spaces $L^p$ for $1<p<\infty$.

\hfill$\Box$

\begin{lem}\label{rottcazztz}
Let $k\geq 2$ be an integer. For every 
$p_3(u)\in {\mathcal P}_3(u)$ such that 
$$\tilde p_3(u)= \prod_{i=1}^3 \partial_x^{\alpha_i} u
\hbox{ with } 0\leq \sum_{i=1}^3 \alpha_i= 2k+1
\hbox{ and } 1\leq \min_{i=1,2,3} 
\alpha_i\leq \max_{i=1,2,3} \alpha_i\leq k $$
we have:
$$\Big |\int p_3(u)dx\Big |\leq C\|u\|_{H^k}^3.$$
\end{lem}

{\bf Proof.} We can assume
$\alpha_1\geq \alpha_2 \geq \alpha_3$ and
also
$p_3(u)=\prod_{i=1}^3 \partial_x^{\alpha_i} u$
(the general case follows in a similar way).
\\
\\
{\em First case: $\alpha_1=\alpha_2=k$}
\\
\\
In this case necessarily 
$\alpha_3=1$ and hence by the H\"older inequality we get
$$\Big |\int p_3(u) dx\Big |\leq \|u\|_{H^k}^2\|\partial_x u\|_{L^\infty}
\leq C \|u\|_{H^k}^3$$
where we have used the Sobolev embedding $H^1\subset L^\infty$.
\\
\\
{\em Second case: $\alpha_2\leq k-1$}
\\
\\
By the H\"older inequality we get
$$\Big |\int p_3(u) dx\Big |\leq \|u\|_{H^k}\|u\|_{H^{k-1}}
\|\partial_x^{\alpha_3}\|_{L^\infty}
$$
and hence by the embedding $H^1\subset L^\infty$
$$...\leq C \|u\|_{H^k}^2\|u\|_{H^{\alpha_3+1}}.$$
The proof follows since
$\alpha_3+1\leq \alpha_2+1\leq k$.

\hfill$\Box$

\begin{lem}\label{quarticstro}
Let $k\geq 1$ and $j\geq 3$ be integers. For every 
$p_j(u)\in {\mathcal P}_j(u)$ such that 
$$\tilde p_j(u)= \prod_{i=1}^j \partial_x^{\alpha_i} u
\hbox{ with } 0\leq \sum_{i=1}^j \alpha_i\leq 2k$$
we have:
\begin{equation}\label{dafar} \Big |\int p_j(u) dx\Big |
\leq C \|u\|_{H^{k}}^j.\end{equation}
\end{lem}

{\bf Proof.} We treat explicitly the case 
$p_j(u)=\tilde p_j(u)$ 
(we specify shortly below how to treat the general case).\\
It is not restrictive to assume that
\begin{equation}\label{eqry}\alpha_1\geq \alpha_2\geq ...\geq \alpha_j.
\end{equation}
By eventually performing integrations by parts we can assume $\alpha_1\leq k$ and by \eqref{eqry}
also $\alpha_2\leq k$. Moreover by the assumption we get
\begin{equation}\label{i>3}\alpha_i<k, \hbox{ } 
\forall i=3,..., j\end{equation}
Hence by the H\"older inequality we get:
\begin{equation}\label{fermi} \Big |\int p_j(u) dx\Big |
\leq \|\partial^{\alpha_1}_x u\|_{L^2} \|\partial_x^{\alpha_2}u\|_{L^2} 
\prod_{i=3}^j \|\partial_x^{\alpha_i}u\|_{L^\infty}
\end{equation}
which due to the embedding $H^1\subset L^\infty$ and \eqref{i>3} implies
\eqref{dafar} 
(if the Hilbert transform $H$ 
is involved in the expression of $p_j(u)$
then we are allowed to remove $H$ at the last step since $\|Hu\|_{H^s}=\|u\|_{H^s}$).

\hfill$\Box$

\begin{lem}\label{lemmacarpiclaim}
Let $n\geq 0$ be an integer and $R>0$, then 
\begin{equation}\label{carpiclaim}
\exists C=C(n,R)>0 \hbox{ s.t. }
\end{equation}
$$\bigcap_{j=0}^{2n} \{u\in H^n||E_{j/2}(u)|< R\}
\subset \{u\in H^n|\|u\|_{H^n}<C
\}.$$
\end{lem}
{\bf Proof.}
We use induction on $n$.
\\
\\
{\em First step: $n=0$}
\\
\\
This is trivial since $E_0(u)=\|u\|_{L^2}^2$.
\\
\\
{\em Second step: $n=1$}
\\
\\
By combining the explicit structure of $E_{1/2}$ (see the introduction)
with the following inequality
$$\|u\|_{L^3}\leq \|u\|_{L^2}^{1/2}\|u\|_{L^6}^{1/2}
\leq C \|u\|_{L^2}^{1/2}\|u\|_{H^{1/2}}^{1/2}$$
we get
\begin{equation}\label{Runmezzo}
|R_{1/2}(u)|\leq 
C \|u\|_{L^2}^{3/2}\|u\|_{H^{1/2}}^{3/2}
\end{equation}
(see the notation in \eqref{strucRmezzi}).
Hence in the region 
$$\{u\in H^1||E_{1/2}(u)|< R, |E_0(u)|<R\}$$
we get
$$\|u\|_{\dot{H}^{1/2}}^2= |E_{1/2}(u) - R_{1/2} (u)|\leq 
R + R^{3/4} \|u\|_{H^{1/2}}^{3/2}$$
which in turn implies the existence of $C>0$ such that
\begin{equation}\label{partcer}
\|u\|_{H^{1/2}}<C, \hbox{ } 
\forall u \in  \{u\in H^1||E_{1/2}(u)|< R, |E_0(u)|<R\}.
\end{equation}
Next, by looking at the explicit structure of $E_1$ (see the introduction)
we get
$$|R_1(u)|\leq C\|u\|_{H^1} \|u\|_{H^{1/2}}^2 + C \|u\|_{H^{1/2}}^4$$
(see the notation \eqref{strucRmezzi})
where we have used the Sobolev embedding
$H^{1/2}\subset L^4$.
Hence by \eqref{partcer} we get a suitable constant $C>0$ such that
$$\|u\|_{\dot{H}^1}^2=|E_1(u) - R_1(u)|\leq R+ C+ C \|u\|_{H^1}$$$$
\hbox{ } \forall u \in  \{u\in H^1||E_1(u)|<R, |E_{1/2}(u)|< R, |E_0(u)|<R\}. $$
In turn, this implies the existence of $C>0$ such that
$$\|u\|_{H^1}< C,
\hbox{ } \forall u \in  \{u\in H^1||E_1(u)|<R, |E_{1/2}(u)|< R, |E_0(u)|<R\}. $$
\\
{\em Third step: $n=2$}
\\
\\
Following the argument of the previous step we get
\begin{equation}\label{snatoro}
\|u\|_{H^1}<C,  \hbox{ } \forall u \in  \{u\in H^2||E_1(u)|<R, |E_{1/2}(u)|< R, |E_0(u)|<R\}
\end{equation}
for a suitable $C>0$.
By combining  the structure of 
$E_{3/2}$ (see the introduction) with \eqref{snatoro} and 
the Sobolev embedding $H^1\subset L^\infty$ we get
$$|R_{3/2}(u)|<C,  \hbox{ } \forall u \in \{u\in H^2||E_1(u)|<R, |E_{1/2}(u)|< R, |E_0(u)|<R\}.$$
As a consequence we deduce
\begin{equation}\label{penjum}\|u\|_{\dot H^{3/2}}^2= |E_{3/2}(u)
-R_{3/2}(u)|\leq R+ C\end{equation}$$ \forall u \in \{u\in H^2||E_{3/2}(u)|<R,
|E_1(u)|<R, |E_{1/2}(u)|< R, |E_0(u)|<R\}.$$
By combining \eqref{penjum} with Lemma \ref{cubichard} and Lemma \ref{quarticstro}
we get
$$\Big |R_2(u) +\frac 54 \int (u_x)^2H u_x \hbox{ } dx\Big |<C$$
$$ \forall u \in \{u\in H^2||E_{3/2}(u)|<R,
|E_1(u)|<R, |E_{1/2}(u)|< R, |E_0(u)|<R\}$$ 
and hence
$$\|u\|_{H^2}^2=E_2(u) + \Big (R_2(u) + \frac 54 \int (u_x)^2H u_x
\hbox{ } dx\Big) 
- \frac 54 \int (u_x)^2H u_x\hbox{ } dx$$
$$\leq R + C + \|u\|_{H^1}^2 \|u\|_{H^2}\leq  R +C +C\|u\|_{H^2}$$
$$\forall u \in \bigcap_{j=0}^4 \{u\in H^2||E_{j/2}(u)|<R\}$$
where we have used the H\"older inequality and the Sobolev embedding
$H^1\subset L^\infty$ to estimate the integral
$\int (u_x)^2H u_x \hbox{ } dx$. 
The proof can be easily concluded.
\\
\\
{\em Fourth step: $n \Rightarrow n+1$ for $n\geq 2$}
\\
\\
Assume the conclusion is proved for $n\geq 1$, then
there exists $C>0$ such that
\begin{equation}\label{induction}
\bigcap_{j=0}^{2(n+1)} \{u\in H^{n+1}||E_{j/2}(u)|< R\}
\subset \bigcap_{j=0}^{2n} \{u\in H^{n+1}||E_{j/2}(u)|< R\} 
\end{equation}$$\subset \{u\in H^{n+1}|\|u\|_{H^{n}}<C\}.
$$
Next we shall use (following \eqref{strucRmezzi}) the notation
\begin{equation}\label{piermer}
E_{n+1/2}(u)= \|u\|_{\dot{H}^{n+1/2}}^2 
+ R_{n+1/2}(u)
\end{equation}
and 
\begin{equation}\label{pietro}
E_{n+1}(u)=\|u\|_{\dot{H}^{n+1}}^2 + R_{n+1}(u)
\end{equation}
(the structure of $R_{k/2}$, described in
\eqref{even} and \eqref{odd}, depending on 
the evenness or the oddness of $k$, will be freely exploited in the sequel).\\
By combining Lemma \ref{quarticstro} (where we choose $k=n$),
with \eqref{induction} we deduce
\begin{equation}\label{pierstro}
|R_{n+1/2}(u) |< C, \hbox{ } 
\forall u\in \bigcap_{j=0}^{2(n+1)} 
\{u\in H^{n+1}||E_{j/2}(u)|< R\}
\end{equation}
for a suitable $C>0$, where we have used the fact that
$R_{n+1/2}(u)$ involves terms of the type $\int p_j(u)dx$
with $j\geq 3$ and $\|p_j(u)\|\leq 2n$
(for a definition of $\|p_j(u)\|$ see \eqref{doublenorm}).
As a consequence of 
\eqref{piermer} and \eqref{pierstro} we get
\begin{equation}\label{induction2}
\bigcap_{j=0}^{2(n+1)} \{u\in H^{n+1}||E_{j/2}(u)|< R\}\subset 
\{u\in H^{n+1}|\|u\|_{H^{n+1/2}}<C\}.
\end{equation}
By combining Lemma \ref{cubichard} (where we choose $m=n$, $\epsilon=1, p=2$),
Lemma \ref{rottcazztz}, Lemma \ref{quarticstro} 
(with $k=n$) and \eqref{induction2} we deduce 
$$|R_{n+1}(u)|< C+ C\|u\|_{H^{n+1}}, \hbox{  } \forall u\in \bigcap_{j=0}^{2(n+1)} \{
u\in H^{n+1}||E_{j/2}(u)|<R\}$$
(where we have used the structure of $R_{n+1}$ given in \eqref{even}).
By combining this estimate with \eqref{pietro} 
we get
$$\|u\|_{\dot H^{n+1}}^2\leq |E_{n+1}(u)| + C + C \|u\|_{H^{n+1}}
\leq R +C +C\|u\|_{H^{n+1}}$$$$
\hbox{ } \forall u\in 
\bigcap_{j=0}^{2(n+1)} \{u\in H^{n+1}||E_{j/2}(u)|< R\}$$
which in turn implies \eqref{carpiclaim} for $n+1$.

\hfill$\Box$

{\bf Proof of Proposition \ref{carpi}}\\
By \eqref{carpiclaim} (where we choose $n=k$) there exists $C>0$ such that  
\begin{equation}\label{sobn-1}
\|u\|_{H^{k}}< C,
\hbox{ } \forall u \in \bigcap_{j=0}^{2k} 
\{u\in H^{k}||E_{j/2}(u)|< R_1\}.
\end{equation}
We also recall the notation (see \eqref{strucRmezzi})
\begin{equation}\label{piermerpri}
E_{k+1/2}(u)= \|u\|_{\dot H^{k+1/2}}^2 
+ R_{k+1/2}(u).
\end{equation}
By combining \eqref{odd} with
\eqref{sobn-1} and 
Lemma \ref{quarticstro}
(recall that we are assuming $k>0$) we get that for every $R$ there exists $C=C(R)$ such that
$$|
R_{k+1/2}(u)|<C, \hbox{ } 
\forall u\in \{u\in H^{k}|\|u\|_{H^{k}}< R\}$$
which is equivalent to
$$|E_{k+1/2}(u) - \|u\|_{\dot H^{k+1/2}}^2|<C, 
\hbox{ } \forall u\in \{u\in H^{k}|\|u\|_{H^{k}}< R\}$$
and hence 
$$|E_{k+1/2}(\pi_N u) -  
\|\pi_N u\|_{\dot H^{k+1/2}}^2
 |<C$$$$ 
\hbox{  } \forall u\in \{u\in H^{k}|\|u\|_{H^{k}}< R\}, 
\hbox{ } N\in \N.$$
By \eqref{sobn-1} we get
$$|
E_{k+1/2}(\pi_N u) - \|\pi_N u\|_{\dot H^{k+1/2}}^2
|<C$$$$
\forall u \in \bigcap_{j=0}^{2k} \{u\in H^{k}||E_{j/2}(u)|< R_1\}, 
\hbox{ } N\in \N$$
that in turn implies 
\eqref{propgian}.

\hfill$\Box$
\section{A linear gaussian bound}
We start with the following general measure theory result which shall be frequently used in the sequel.
\begin{prop}\label{expon2}
Let $F:(\Omega, {\mathcal A}, p)\rightarrow {\mathbf C}$  be measurable and $C, \alpha>0$ be such that
\begin{equation}\label{tchba}\|F\|_{L^q}\leq C q^\alpha, \hbox{ } \forall q\in [1, \infty).\end{equation}
Then $$p\{\omega\in \Omega||F(\omega)|>\lambda\}\leq e^{-\frac{\alpha}e \left(
\frac \lambda C\right)^\frac 1\alpha}, \hbox{ } \forall \lambda>0.$$ 
\end{prop}
{\bf Proof.}
By combining the Tchebychev inequality 
with \eqref{tchba} we get:
$$p\{\omega\in \Omega||F(\omega)|>\lambda\}\leq \frac{\|F\|_q^q}{\lambda^q} 
\leq C^q \Big ( \frac {q^\alpha} \lambda\Big)^q.$$
We conclude by choosing
$q =\Big ( \frac{\lambda}{C}\Big )^\frac 1\alpha e^{-1}.$

\hfill$\Box$

Next we present, as an application of the previous result, a linear gaussian bound which will be used in next sections.
\begin{prop}\label{sobo2}
For every integer $m\geq 0$, $0<\epsilon<\frac 12$ and $1\leq p <\infty$
there exists $C=C(m,\epsilon, p)>0$ such that
$$p \{\omega\in \Omega| 
\|\varphi_{m+1}(\omega)
\|_{W^{\epsilon, p}_x}
>\lambda \} \leq  C e^{-\frac{\lambda^2}{C}},
\hbox{ } \forall \lambda>0$$
where $\varphi_{m+1}(\omega)$ is the random vector in \eqref{randomized}
for $k=2(m+1)$.
\end{prop}
{\bf Proof.}
It is sufficient to prove that
$$p\Big\{\omega\in \Omega| \Big \|
\sum_{k>0} \frac{\varphi_k(\omega)}
{k^{m+1-\epsilon}} e^{ikx}
\Big \|_{L^p_x}
>\lambda \Big \} \leq C e^{-\frac{\lambda^2}C}.$$
For every 
fixed $x\in (0, 2\pi)$ 
the random variable 
$$\sum_{k>0} 
\frac{\varphi_k(\omega)}{k^{m+1-\epsilon}} e^{ikx}$$
is Gaussian and its distribution function is
$\frac 1{\pi K}e^{-\frac{ |z|^2}{K}} dz$
where
$K=\sum_{k>0} \frac 1{k^{2(m+1-\epsilon)}}.$\\
As a consequence we get the following estimate:
$$\Big \|\sum_{k>0} 
\frac{\varphi_k(\omega)}{k^{m+1-\epsilon}} 
e^{ikx}\Big \|_{L^q_\omega}^q
= \frac 1{\pi K}\int_{\C} |z|^q e^{-\frac{|z|^2}K} dz$$ 
$$=2K^\frac q2 \int_0^\infty e^{-s^2} s^{q+1} ds\leq C
K^\frac q2 \left( \frac q2\right )^\frac q2,
\hbox{ } \forall x\in (0,2\pi)$$
for a suitable $C>0$ (the last inequality 
can be proved by 
integration by parts).
In particular we get:
$$\Big \|\sum_{k>0} \frac{\varphi_k(\omega)}
{k^{m+1-\epsilon}} e^{ikx}
\Big \|_{L^q_\omega}
\leq C \sqrt q, \hbox{ } \forall x\in (0, 2\pi)$$
and
hence
$$\Big \|\sum_{k>0} \frac{\varphi_k(\omega)}
{k^{m+1-\epsilon}} e^{ikx}
\Big \|_{L^p_x L^q_\omega}
\leq C \sqrt q.$$
Due to the inequality 
$\|.\|_{L^q_\omega L^p_x}\leq \|.\|_{L^p_xL^q_\omega}$
for every $q\geq p$ we get:
$$\Big \|
\sum_{k>0} \frac{\varphi_k
(\omega)}{k^{m+1-\epsilon}} e^{ikx}\Big \|_{L^q_\omega L^p_x}
\leq C \sqrt q, \hbox{ } \forall q\geq p$$
(since $p(\Omega)=1$ it is easy to deduce that
the estimate above is true for every 
$q\geq 1$ eventually with a
new constant $C$).
Hence we can conclude by using Proposition~\ref{expon2}.

\hfill$\Box$

\section{Multilinear gaussian bounds}
For any $p(u)\in \cup_{j=1}^\infty{\mathcal P}_j(u)$
(see Section \ref{conserlawsbo}) and for any $N\in \N$ we introduce the functions
$$f^p(v)=\int p(v) dx \hbox{ and }
f^p_N(v)= \int p(\pi_N v) dx.$$
We also recall that the Sobolev spaces $H^{m+1/2-\epsilon}$
are a support for the Gaussian measure
$d\mu_{m+1}$ for every $\epsilon>0$. This fact will be used 
without any further comment in the sequel.
The main results of this section are the following propositions.
\begin{prop}\label{wiener} Let $m\geq 0$ and 
$p_3(u)\in {\mathcal P}_3(u)$ 
be such that 
$\tilde p_3(u)= u \partial_x^m u 
\partial_x^{m+1}u.$
Then there exists $C>0$ such that 
\begin{equation}\label{2pwie}
\|f_{N}^{p_3}(u)- f_{M}^{p_3}(u)\|_{L^p(d\mu_{m+1})}
\leq C \frac {p^{3/2}}{\sqrt {\min\{M,N\}}},  
\hbox{ } \forall M,N\in \N, p\geq 2.
\end{equation}
In particular 
\begin{equation}\label{third}
\exists C>0 \hbox{ s.t. } 
\mu_{m+1}(A_{M,N}^{p_3,\lambda})\leq  
e^{-\frac1C (\lambda \sqrt {\min \{N,M\}})^{2/3}}, \hbox{ }
\forall M,N\in \N, \lambda>0
\end{equation}
where
\begin{equation}\label{Amn}
A_{M,N}^{p_3,\lambda}=\{u\in H^{m+1/2 -\epsilon}
||f_{N}^{p_3}(u)- f_{M}^{p_3}(u)|>\lambda\}.
\end{equation}
\end{prop}
\begin{prop}\label{wiener2} Let $m\geq 0$ be fixed. There exists $C>0$ such that
\begin{equation}\label{2pwie2}
\|h_{N}(u)- h_{M}(u)\|_{L^p(d\mu_{m+1})}
\leq C \frac {p}{\sqrt {\min\{M,N\}}},  
\hbox{ } \forall M,N\in \N, p\geq 2
\end{equation}
where $h_K(u)= \|\pi_K u\|_{\dot H^{m+1/2}}^2 -\alpha_K$ for any $K\in \N$.
In particular 
\begin{equation}\label{fourth}
\exists C>0 \hbox{ s.t. } \mu_{m+1}(B_{M,N}^\lambda)\leq  
e^{-\frac1C (\lambda \sqrt {\min \{N,M\}})},
\hbox{ }
\forall M,N\in \N, \lambda>0
\end{equation}
where
\begin{equation}\label{Bmn}
B_{M,N}^\lambda=\{u\in H^{m+1/2-\epsilon}
| |h_{N}(u) - h_{M}(u)|>\lambda
\}.\end{equation}
\end{prop}
We need some preliminary lemmas. The first one concerns the orthogonality of the functions 
$\{\varphi_i(\omega) \varphi_j(\omega) \varphi_k(\omega)\}_{i,j,k\in \Z}$
(where $\varphi_n(\omega)$ are the Gaussian functions
that appear in \eqref{randomized})
provided that
$(i,j,k)\in {\mathcal A}$, where
\begin{equation}\label{rabo}
{\mathcal A}=\{(i,j,k)\in \Z\setminus\{0\} | i+j+k=0\}.
\end{equation}
\begin{lem}\label{orth}
Let
$(k_1, k_2, k_3), (j_1, j_2, j_3)\in \mathcal A$
be such that
\begin{equation}\label{hplesa}
\{k_1, k_2, k_3\}
\neq \{j_1, j_2, j_3\}.\end{equation}
Then
$$\int_{\Omega} \varphi_{k_1}\varphi_{k_2}\varphi_{k_3}
\overline{\varphi_{j_1}\varphi_{j_2}\varphi_{j_3}} d\omega=0.$$ 
\end{lem}
{\bf Proof.}  We split the proof in two cases 
(which in turn are splited in several subcases).\\
\\
{\bf First case}: $\exists i\in \{1,2,3\}$ {\em s.t.} 
$k_i \notin \{j_1, j_2, j_3\}$
\\
\\
We can assume 
\begin{equation}\label{hpfirstcase}
k_1 \notin \{j_1, j_2, j_3\}.
\end{equation}
Next we consider four subcases:
\\

{\em First subcase:}
\begin{equation}\label{fcpe}k_1 \notin \{k_2, k_3\} \hbox{ and } 
-k_1 \notin \{j_1, j_2, j_3\}\end{equation}
Notice that by definition of $\mathcal A$ necessarily
\begin{equation}\label{hfc}
-k_1\notin  \{k_2, k_3\}.\end{equation} Hence
by combining \eqref{hpfirstcase}, \eqref{fcpe}, \eqref{hfc} and 
the independence assumption on $\{\varphi_n(\omega)\}_{n>0}$ we get:
$$0=\int \varphi_{k_1}\int \varphi_{k_2}\varphi_{k_3}
\overline{\varphi_{j_1}\varphi_{j_2}\varphi_{j_3}} d\omega=
\int \varphi_{k_1} \varphi_{k_2}\varphi_{k_3}
\overline{\varphi_{j_1}\varphi_{j_2}\varphi_{j_3}} d\omega.$$

{\em Second subcase}: \begin{equation}\label{secsub0}
k_1 \in \{k_2, k_3\}
\hbox{ and } 
-k_1 \notin \{j_1, j_2, j_3\}\end{equation}
It is not restrictive to assume 
\begin{equation}\label{ident9}k_1=k_2\end{equation} and
hence by the definition of  $\mathcal A$
\begin{equation}\label{secsub}
k_3\neq \pm k_1.
\end{equation}
Hence by combining \eqref{hpfirstcase}, 
\eqref{secsub0}, \eqref{ident9},
\eqref{secsub}
with the independence assumption on $\{\varphi_n(\omega)\}_{n>0}$ we get
$$0=\int \varphi_{k_1}^2 d\omega \int \varphi_{k_3}
\overline{\varphi_{j_1}\varphi_{j_2}\varphi_{j_3}} d\omega=
\int \varphi_{k_1} \varphi_{k_2}\varphi_{k_3}
\overline{\varphi_{j_1}\varphi_{j_2}\varphi_{j_3}} d\omega.$$

{\em Third subcase}: 
\begin{equation}\label{thrthi}
k_1 \notin \{k_2, k_3\}
\hbox{ and } 
-k_1 \in \{j_1, j_2, j_3\}
\end{equation}
By definition of $\mathcal A$ we also deduce
\begin{equation}\label{caz}
-k_1 \notin \{k_2, k_3\}.\end{equation} By \eqref{thrthi}
we can assume for simplicity
either 
\begin{equation}\label{1}-k_1=j_1\notin \{j_2, j_3\}\end{equation}
or 
\begin{equation}\label{2}-k_1=j_1= j_2.\end{equation}
In the case when \eqref{1} occurs, we can also assume
by the definition of $\mathcal A$ that 
\begin{equation}\label{lululu}
k_1=-j_1\notin \{j_2, j_3\}.
\end{equation}
By combining \eqref{hpfirstcase}, \eqref{thrthi}, \eqref{caz}, 
\eqref{1}, \eqref{lululu} with the independence assumption
on $\{\varphi_n(\omega)\}_{n>0}$ 
we get$$0=\int \varphi_{k_1}^2 d\omega \int \varphi_{k_2}
{\varphi_{k_3}\overline {\varphi_{j_2}\varphi_{j_3}}} d\omega=
\int \varphi_{k_1} \varphi_{k_2}\varphi_{k_3}
\overline{\varphi_{j_1}\varphi_{j_2}\varphi_{j_3}} d\omega$$
(where we have used also \eqref{hpfirstcase});\\
in the case when \eqref{2} occurs, by using the definition of 
$\mathcal A$ we get
\begin{equation}\label{lul}
\pm k_1\neq j_3. 
\end{equation}
Hence
by combining 
\eqref{thrthi}, \eqref{caz}, \eqref{2}, \eqref{lul} we deduce
$$0=\int \varphi_{k_1}^3 d\omega \int \varphi_{k_2}
{\varphi_{k_3}\overline {\varphi_{j_3}}} d\omega  =
\int \varphi_{k_1} \varphi_{k_2}\varphi_{k_3}
\overline{\varphi_{j_1}\varphi_{j_2}\varphi_{j_3}} d\omega$$

{\em Fourth subcase}: \begin{equation}\label{471}
k_1 \in \{k_2, k_3\}
\hbox{ and }
-k_1 \in \{j_1, j_2, j_3\}\end{equation}
We can assume
\begin{equation}\label{789}
k_1=k_2
\end{equation}
and by the definition of $\mathcal A$
also
\begin{equation}\label{8910}k_3\neq \pm k_1.\end{equation}
Moreover, we can assume that
either 
\begin{equation}\label{1bis}-k_1=j_1\notin 
\{j_2, j_3\}\end{equation}
or 
\begin{equation}\label{2bis}-k_1=j_1= j_2.\end{equation}
In the case when \eqref{1bis} occurs we can also assume
by the definition of $\mathcal A$ that 
\begin{equation}\label{patval}
k_1=-j_1\notin \{j_2, j_3\}.\end{equation} 
Hence by combining 
\eqref{789}, \eqref{8910}, \eqref{1bis},\eqref{patval}
we get $$0=\int \varphi_{k_1}^3 d\omega \int
{\varphi_{k_3}\overline {\varphi_{j_2}\varphi_{j_3}}} d\omega=
\int \varphi_{k_1} \varphi_{k_2}\varphi_{k_3}
\overline{\varphi_{j_1}\varphi_{j_2}\varphi_{j_3}} d\omega;$$
in the case when \eqref{2bis} occurs we can deduce by the definition
of $\mathcal A$ that 
\begin{equation}\label{vadi}
\pm k_1\neq j_3.                      
\end{equation}
Hence by \eqref{789}, \eqref{8910}, \eqref{2bis}, \eqref{vadi}
we get
$$0=\int \varphi_{k_1}^4 d\omega \int
{\varphi_{k_3}\overline {\varphi_{j_3}}} d\omega  =
\int \varphi_{k_1} \varphi_{k_2}\varphi_{k_3}
\overline{\varphi_{j_1}\varphi_{j_2}\varphi_{j_3}} d\omega.$$
\\
{\bf Second case:}
\begin{equation}\label{secosuc} 
k_i \in\{j_1, j_2, j_3\}
\hbox{ } \forall i=1,2,3\end{equation}
\\
\\
Next we consider two subcases:
\\

{\em First subcase}: 
\begin{equation}\label{fevad}k_1\neq k_2, k_1\neq k_3,
k_2\neq k_3\end{equation}
\\
By combining \eqref{secosuc} and \eqref{fevad} 
it is easy to deduce that
$$\{j_1, j_2, j_3\}=\{k_1, k_2, k_3\}$$
which is in contradiction with \eqref{hplesa}.
\\

{\em Second subcase:}
\begin{equation}\label{secosu} 
\exists n, m\in \{1,2,3\}
 \hbox{
 s.t. } n\neq m, k_n=k_m\end{equation}
\\
We can assume 
\begin{equation}\label{secoug}
k_1=k_2,\end{equation} then by the definition of $\mathcal A$
we deduce that
\begin{equation}\label{secugu}
k_3=-2k_1.\end{equation}
On the other hand by \eqref{secosuc} 
$k_1, -2 k_1\in \{j_1, j_2, j_3\}$.
Since by the definition of $\mathcal A$ we have $\sum_{i=1}^3 j_i=0$
we conclude that necessarily 
\begin{equation}\label{usgualiesplicit}
\{j_1, j_2, j_3\}=\{k_1, k_1, -2k_1\}.
\end{equation}
On the other hand by \eqref{secoug}, \eqref{secugu}
we get $$\{k_1, k_2, k_3\}=\{k_1, k_1, -2k_1\}$$
which in conjunction with \eqref{usgualiesplicit} gives
$\{k_1, k_2, k_3\}=\{j_1, j_2, j_3\}$. 
Hence we get a contradiction with 
the hypothesis \eqref{hplesa}.
\hfill$\Box$
\begin{lem}\label{orthogonal}
Let $m\geq 0$ be an integer and $p_3(u)\in {\mathcal P}_3(u)$
such that $\tilde p_3(u)=u\partial_x^mu\partial_x^{m+1}u.$
Then there exists $C>0$ such that
$$\\|f_{N}^{p_3}(u)- f_{M}^{p_3}(u)\|_{L^2(d\mu_{m+1})}\leq 
\frac C{\sqrt {\min\{N,M\}}}, \hbox{ } \forall N,M\in \N.
$$
\end{lem} 
{\bf Proof.} We assume for simplicity $p_3(u)=u\partial_x^m u \partial_x^{m+1}u$
(the general case can be treated in a similar way).
Next we  assume $N> M$ and we shall use the parametrization \eqref{randomized}
with $k=2(m+1)$ to describe our probability space.
Hence we get the following representation
$$f^{p_3}_{N}(\varphi(\omega))-f_{M}^{p_3}(\varphi(\omega))=
\sum_{(i,j,k)\in {\mathcal A}_{M}^N} 
\frac{1}{|i|^{m+1}} \frac 1{|j|} \varphi_i(\omega) 
\varphi_j(\omega) \varphi_k(\omega)$$
where
\begin{equation}\label{Astorto}{\mathcal A}_{M}^N=\{(i,j,k)\in {\mathcal A}
| |i|, |j|, |k|\leq N
\hbox{ and } \max\{|i|, |j|, |k|\}>M\}\end{equation}
and $\mathcal A$ is defined as in \eqref{rabo}.
By Lemma \ref{orth} we get:
$$\|f_{N}^{p_3}(\varphi(\omega))-f_{M}^{p_3}
(\varphi(\omega))\|_{L^2_\omega}^2$$$$\leq
\sum_{(i,j,k)\in \tilde {\mathcal A}_M^N} 
\Big [\frac{1}{|i|^{m+1}}\Big (
\frac 1{|j|}+\frac{1}{|k|}\Big ) 
+\frac{1}{|j|^{m+1}}
\Big (\frac 1{|i|}+\frac{1}{|k|}\Big ) 
+\frac{1}{|k|^{m+1}}\Big (\frac 1{|i|}
+\frac{1}{|j|}\Big )
\Big]^2$$ where
$$\tilde {\mathcal A}_{M}^N=\{(i,j,k)\in 
{\mathcal A}_{M}^N|i\leq j\leq k\}.$$
Next notice that the following 
elementary property holds
$$\tilde {\mathcal A}_{M}^N\subset \{(i,j,k)
\in {\mathcal A}| \hbox{Card}\{|i|,|j|,|k|\}\in 
[M/2, N] \}\geq 2\}$$
and hence we easily get:
$$\sum_{(i,j,k)\in \tilde {\mathcal A}_M^N} 
\Big [\frac{1}{|i|^{m+1}}\Big (\frac 1{|j|}
+\frac{1}{|k|}\Big ) 
+\frac{1}{|j|^{m+1}}
\Big (\frac 1{|i|}+\frac{1}{|k|}\Big ) 
+\frac{1}{|k|^{m+1}}\Big (\frac 1{|i|}
+\frac{1}{|j|}\Big )
\Big ]^2$$
$$\leq C \sum_{\substack{(l,n)\in \N\times \N
\\n\geq \frac M2}} 
\frac 1{l^{2(m+1)}n^2}\leq \frac CM.$$

\hfill$\Box$

In next lemma the functions $h_K(u)$ are the ones defined in Proposition
\ref{wiener2}.
\begin{lem}\label{orthogonal2}
Let $m\geq 0$ be an integer. Then there exists $C>0$ such that
$$\|h_{N}(u)- h_{M}(u)
\|_{L^2(d\mu_{m+1})}\leq \frac C{\sqrt {\min\{N,M\}}}, \hbox{ } 
\forall N,M\in \N.$$
\end{lem}
{\bf Proof.} Notice that
$\|\varphi(\omega)\|_{\dot H^{m+ 1/2}}^2
= \sum_{n\in \Z\setminus \{0\}}
\frac1{|n|} |\varphi_n(\omega)|^2
$ where $\varphi(\omega)$ is defined as in \eqref{randomized}
for $m=2(k+1)$. Hence the proof 
follows as in \cite{tz} (see Lemma 4.7).

\hfill$\Box$

{\bf Proof of Proposition \ref{wiener}} 
In Lemma~\ref{orthogonal} we have proved 
\eqref{2pwie} for $p=2$.
The case $p>2$ follows by combining 
the estimate for $p=2$  
with the Wiener Chaos
in the same spirit as the paper \cite{tz} (see the proof of Lemma ~4.3 in \cite{tz}).
The estimate \eqref{third} follows by 
\eqref{2pwie} in conjunction
with Proposition \ref{expon2}.

\hfill$\Box$
 
 We refer to \cite{LT} for a background on the estimates for the Wiener Chaos.
 \\
 
{\bf Proof of Proposition \ref{wiener2}} 
By combining Lemma \ref{orthogonal2} 
with the Wiener Chaos in the spirit 
of \cite{tz} we get \eqref{2pwie2}
for any $p\geq 2$ (see the proof of Lemma~4.8 in \cite{tz}).
Finally \eqref{fourth} follows by combining \eqref{2pwie2} 
with Proposition~\ref{expon2}.

\hfill$\Box$

Arguing as in the proof of 
Proposition \ref{wiener} and \ref{wiener2}
we can prove the following result (that will be useful in the 
sequel to prove Theorem~\ref{main}
in the special case $k=2$).
\begin{prop}\label{new}
There exists $C>0$ such that
\begin{equation}\label{firstnew}
\mu_1\{u\in H^{1/2 -\epsilon}| |E_{1/2}(\pi_N u) - 
\alpha_N -E_{1/2}(\pi_M u) +\alpha_M|>\lambda\}
\end{equation}$$\leq e^{-\frac 1C (\lambda\sqrt{\min\{M,N\}})^{2/3}}$$
and 
\begin{equation}\label{firstnew2}
\mu_1\{u\in H^{1/2 -\epsilon}| \|\pi_N u - 
\pi_M u\|^2_{L^2}>\lambda\}
\leq e^{-\frac 1C (\lambda\sqrt{\min\{M,N\}})}
\end{equation}
$$\hbox{ } \forall M, N\in \N, \lambda>0.$$
\end{prop}
{\bf Proof.} The proof of \eqref{firstnew2}
follows the same argument as the proof of \eqref{third}
and \eqref{fourth} (i.e. it follows by combining Lemma \ref{orth}, with the Wiener Chaos
and Proposition \ref{expon2}).
By a similar argument we can prove
$$\mu_1\{u\in H^{1/2 -\epsilon}| |R_{1/2}(\pi_N u) - R_{1/2}(\pi_Mu)|>\lambda\}
\leq e^{-\frac 1C (\lambda\sqrt{\min\{M,N\}})^{2/3}}.$$
By combining this estimate with \eqref{fourth} (for $m=0$)
we get \eqref{firstnew}.
\hfill$\Box$
\section{Proof of Theorem \ref{main} for $k=2(m+1)$,  $m\geq 0$}
Along this section, when it is not better specified, we shall assume that $m\geq 0$ is a given integer.
We recall the following notations to describe the energies preserved by the Benjamin-Ono flow:
$$E_{m+1}(u)=\|u\|_{\dot H^{m+1}}^2
+R_{m+1}(u);$$
$$E_{m+ 1/2}(u)=\|u\|_{\dot H^{m+1/2}}^2
+R_{m+1/2}(u).$$
We also introduce the following functions
$$f_{N}:H^{m+1/2-\epsilon}\ni u\rightarrow  R_{m+1} (\pi_N u);$$
$$g_{N}:  H^{m+1/2 -\epsilon}\ni u\rightarrow 
E_{m+1/2} (\pi_N u) 
-\alpha_N$$
(recall that $H^{m+1/2 -\epsilon}$ is of full measure for $\mu_{m+1}$).
Notice that we can write the identity
$$g_{N}(u) - h_{N}(u)=
R_{m+1/2} (\pi_N u)$$
where $h_{N}(u)$ is defined as in Proposition \ref{wiener2}.
For every $p_j(u)\in {\mathcal P}_j(u)$
we introduce
\begin{equation}
f_{N}^{p_j}: H^{m+ 1/2-\epsilon} \ni u\rightarrow 
 \int p_j(\pi_N u) dx\in \R.
\end{equation}
Next we split the proof of Theorem \ref{main}
(in the case $k=2(m+1)$) 
in several propositions.

\begin{prop}\label{proprImain}
Let $m\geq 0$ be an integer and $\psi\in C_c(\R)$ be given. 
Then there exist two
functions $\bar h(u), \bar f(u)$
measurable with respect to $\mu_{m+1}$
such that: 
\begin{equation}\label{finiti}
|\bar h (u)|, |\bar f(u)|<\infty, \hbox{ a.e. (w.r.t. $\mu_{m+1}$) } u\in H^{m+ 1/2 - \epsilon};
\end{equation}
\begin{equation}\label{infiniti}\prod_{j=0}^{2m} \psi(E_{j/2}(\pi_N u))
\psi( E_{m+1/2}(\pi_N u) - \alpha_N )e^{-R_{m+1}(\pi_N u)}
\end{equation}$$ \hbox{ converges in measure to }$$
$$\prod_{j=0}^{2m} \psi(E_{j/2}(u))
\psi(\bar h(u) + R_{m+1/2}(u)) 
e^{-\bar f(u)}.$$
Moreover
\begin{equation}\label{lastbnotleas}
|E_{j/2}(u)|, |R_{m+1/2}(u)|<\infty,
\hbox{ a.e. (w.r.t. $\mu_{m+1}$) }  u\in   H^{m+ 1/2 - \epsilon}. 
\end{equation}
\end{prop}
The proof of \eqref{lastbnotleas} 
follows by \eqref{inpart} and \eqref{ER} in Lemma~\ref{Rbar}. Hence 
Proposition~\ref{proprImain} follows by
Lemmas~\ref{Rbar}, \ref{Rbar2} in conjunction with the following
proposition.
\begin{prop}\label{barfh}
Let $m\geq 0$ be an integer and $\psi\in C_c(\R)$ be given. There exist 
$\bar f(u), \bar h(u)$ measurable functions 
with respect to $\mu_{m+1}$ such that:
\begin{equation}\label{sofbenmon}|\bar h(u)|, |\bar f(u)|<\infty,
\hbox{ a.e. (w.r.t. $\mu_{m+1}$) } u \in H^{m+ 1/2 - \epsilon};
\end{equation} 
\begin{equation}\label{barf}
R_{m+1}(\pi_N u) \hbox{ converges in measure w.r.t. $\mu_{m+1}$ to } 
\bar f(u);
\end{equation}
\begin{equation}\label{barh}
\lim_{N\rightarrow \infty}
\|\psi(g_{N}(u))- \psi(\bar h(u) + 
R_{m+1/2}(u))\|_{L^q(d\mu_{m+1})}=0, \hbox{ } \forall q\in [1, \infty).
\end{equation}
Moreover we have
\begin{equation}\label{truncaturelimitcaz}
\lim_{N\rightarrow \infty} 
\|\psi(E_{j/2} (\pi_N u))
- \psi(E_{j/2}(u))\|_{L^q(d\mu_{m+1})}=0,
\hbox{ } \forall q\in [1, \infty), j=0,...,2m.
\end{equation}
\end{prop}
First we prove the following lemma.
\begin{lem}\label{Rbar}
Let $m\geq 1$ be an integer 
then the following limits exist:
\begin{equation}\label{noninpart}\lim_{N\rightarrow \infty}
f_{N}^{p_j}(u)= \int p_j(u) dx
\in \R, \hbox{ a.e. (w.r.t. $\mu_{m+1}$) } u\in H^{m+ 1/2 -\epsilon}
\end{equation}
provided that \begin{equation}\label{0471}
j\geq 3, \hbox{ } \tilde p_j(u)= \prod_{i=1}^j \partial_x^{\alpha_i} u
\hbox{ with } 0\leq \sum_{i=1}^j \alpha_i\leq 2m.
\end{equation}
or
\begin{equation}\label{0477}
j=3, \hbox{ } \tilde p_3(u)= \prod_{i=1}^3 \partial_x^{\alpha_i} u
\hbox{ with } 0\leq \sum_{i=1}^3 \alpha_i= 2m+1
\hbox{ and } \min_{i=1,2,3} \alpha_i\geq 1.
\end{equation}
In particular
\begin{equation}\label{inpart}
\lim_{N\rightarrow \infty} R_{m+1/2} (\pi_N u)
=R_{m+1/2}(u)\in \R\end{equation}$$ \hbox{ a.e. (w.r.t. $\mu_{m+1}$) } 
u\in H^{m+ 1/2 - \epsilon}.$$
Moreover we have
\begin{equation}\label{ER}
|E_{j/2}(u)|<\infty, \hbox{ } \forall u\in H^{m+ 1/2 -\epsilon}, j=0,..., 2m.  
\end{equation}
\end{lem}
{\bf Proof.} 
We assume  for simplicity $p_j(u)=\tilde p_j(u)$
(the general case can be treated by a similar argument).\\ 
The proof of \eqref{noninpart} (under the hypothesis \eqref{0471})
follows by  
Lemma \ref{quarticstro}.\\
Concerning the proof of
\eqref{noninpart}, under the assumption \eqref{0477}, we
notice that 
by integration by parts we can assume 
$$1\leq \alpha_1\leq \alpha_2\leq \alpha_3 \leq m.$$
Hence we get
$$\Big |\int p_3(u) dx\Big |\leq C \|u\|_{H^m}^2\|\partial_x^{\alpha_1} u\|_{L^\infty}
\leq C \|u\|_{H^m}^2\|u\|_{W^{m+\epsilon,p}}$$
where we have used the Sobolev embedding
$W^{\epsilon, p} \subset L^\infty$ provided that $\epsilon p>1$.
On the other hand by a suitable version of Proposition \ref{sobo2}
(where we replace $\|.\|_{W^{\epsilon,p}}
$by $\|.\|_{W^{m+\epsilon, p}}$) we get 
$u\in W^{m+\epsilon, p}, \hbox{ a.e. ( w.r.t. $\mu_{m+1}$) } u\in H^{m+1/2 - \epsilon}$
and hence \eqref{noninpart} follows.\\
The proof of \eqref{inpart} follows by combining the structure
of $E_{m+1/2}$ (see \eqref{odd})
with \eqref{noninpart} (under the assumption \eqref{0471}).
The proof of \eqref{ER} follows by a similar argument.

\hfill$\Box$

The next result is a suitable version of the previous lemma in the case $m=0$.
\begin{lem}\label{Rbar2}
The following limits exist:
\begin{equation}\label{4545}
\lim_{N\rightarrow \infty}
R_{1/2} (\pi_N u)=R_{1/2}(u)\in \R, \hbox{ } \forall u \in H^{1/2 -\epsilon};
\end{equation}
\begin{equation}\label{4444}
\lim_{N\rightarrow \infty} \int (\pi_N u)^4 dx
=\int u^4 dx, \hbox{ } \forall u \in H^{1/2-\epsilon}.
\end{equation}
\end{lem}

{\bf Proof.} 
By looking at the explicit structure
on $E_{1/2}$ (see the introduction)
we get $R_{1/2}(u)=\frac 13 \int u^3 dx$.
On the other hand by the Sobolev embedding
$H^{1/2-\epsilon}\subset L^3$
we get $u\in L^3$, and hence \eqref{4545} follows. By a similar argument we deduce
\eqref{4444}.

\hfill$\Box$

{\bf Proof of Proposition \ref{barfh}}
By Proposition~\ref{wiener}
and \ref{wiener2} there exist two functions $\tilde f, 
\bar h\in \cap_{q=1}^\infty L^q(d\mu_{m+1})$ 
such that:
\begin{equation}\label{Ibache}
\lim_{N\rightarrow \infty} \|f^{p_3}_{N}(u)-\tilde f(u)
\|_{L^q(d\mu_{m+1})}=0
\hbox{ (provided that }
\tilde p_3(u)=u\partial_x^m u\partial_x^{m+1}u);
\end{equation}
\begin{equation}\label{IIbache}
\lim_{N\rightarrow \infty} \|h_{N}(u)-\bar h(u)\|_{L^q(d\mu_{m+1})}=0.
\end{equation}
{\bf Proof of \eqref{barf}}
\\
\\
If $m=0$ then it follows by \eqref{4444},
\eqref{Ibache}
and by looking at the explicit structure of $E_1$ (see the introduction).\\
If $m\geq 1$ then it follows by combining
\eqref{noninpart} (under both assumptions
\eqref{0471} and \eqref{0477}), \eqref{Ibache} and
the algebraic structure of $R_{m+1}(u)$ (see \eqref{even}).
\\
\\
{\bf Proof of \eqref{barh}}
\\
\\
It is sufficient to prove that for every 
sequence $N_k$ there is a subsequence
$N_{k_h}$ such that 
\begin{equation}\label{subsequence}
\lim_{h\rightarrow \infty}
\|\psi(g_{N_{k_h}}(u))- \psi(\bar h(u) + 
R_{m+1/2}(u))\|_{L^q(d\mu_{m+1})}=0.
\end{equation}
Notice that by combining \eqref{inpart} (when $m\geq 1$)
and \eqref{4545} (when $m=0$) with \eqref{IIbache} we get a subsequence $N_k$
such that 
$$\lim_{k\rightarrow \infty} \Big(h_{N_k}(u)
+ R_{m+1/2} (\pi_{N_k} u)\Big)
=\bar h(u) +R_{m+1/2}(u)$$$$ \hbox{ a.e. (w.r.t. $\mu_{m+1}$) } 
u \in H^{m+1/2 -\epsilon}.$$
Since 
$\sup_{u\in H^{m+ 1/2 -\epsilon}}|\psi(g_{N_k}(u))|\leq \sup \psi$
and $d\mu_{m+1}(H^{m+ 1/2 -\epsilon})=1$
we can apply the dominated convergence theorem to get
\eqref{subsequence}.
\\
\\
{\bf Proof of \eqref{truncaturelimitcaz}}
\\
\\
If $m\geq 1$ then we combine 
\eqref{noninpart} (under the assumption \eqref{0471})
with \eqref{even} and \eqref{odd} in order to get
$E_{j/2}(\pi_N u)\rightarrow
E_{j/2}(u), \hbox{ a.e. (w.r.t. $\mu_{m+1}$) } u\in H^{m+ 1/2 - \epsilon}$
for $j=0,.., 2m$.
Hence the proof for $m\geq 1$ can be concluded as in \eqref{barh}.\\
The case $m=0$ is simpler
since we have $\mu_1(L^2)=1$ and hence
$\|\pi_N u\|_{L^2}\rightarrow \|u\|_{L^2},
\hbox{ a.e. (w.r.t. $\mu_1$) } u\in L^2.
$ 
The proof follows as above.

\hfill$\Box$

The next proposition allows us to deduce that
the limit functions constructed in \eqref{infiniti}
belong to $L^q(d\mu_{m+1})$.
\begin{prop}\label{proprImain45}
Let $m\geq 0$ and $\psi\in C_c(\R)$ be given. 
For every $q\in [1,\infty)$ we have
\begin{equation}\label{incorporame}
\sup_N \Big \|\prod_{j=0}^{2m} \psi(E_{j/2}(\pi_N u))
\psi(E_{m+1/2}(\pi_N u) -\alpha_N )e^{-R_{m+1}(\pi_N u)}
\Big \|_{L^q(d \mu_{m+1})}<\infty.
\end{equation}
\end{prop}
\begin{lem}\label{degreg}
Let $m\geq 0$ be an integer and $p_3(u)\in {\mathcal P}_3(u)$ such that
$\tilde p_3(u)=u\partial_x^mu \partial_x^{m+1}u.$
For every $R>0$ there exists $C=C(R)>0$ such that
\begin{equation}\label{first}
\mu_{m+1}\{u\in H^{m+1/2 -\epsilon} | 
|f_{N}^{p_3}(u)|>\lambda, |h_{N}(u)|<R\}
\leq C e^{-\frac{\lambda^2}{C \alpha_N^2}}\end{equation}
$$\forall N\in \N, \lambda>0.
$$
\end{lem}
{\bf Proof.} 
We fix $0<\epsilon<\frac 12$ and $1< p<\infty$ 
such that $\epsilon p>1$. Then by Lemma \ref{cubichard} we get
$$|f_{N}^{p_3}(u)|\leq C \|\pi_N
u\|_{H^{m+1/2}}^2
\|\pi_N u\|_{W^{\epsilon, p}}
$$
and hence
$$|f_{N}^{p_3}(u)|\leq C (\alpha_N+R)
\|\pi_N u\|_{W^{\epsilon,p}}
$$
$$\forall u\in \{u\in H^{m+1/2 -\epsilon}|
|h_{N}(u)|<R\}.$$
The proof follows by Proposition \ref{sobo2}
(in fact notice that the same proof of Proposition~\ref{sobo2} works in the case when  
the vector $\varphi(\omega)$ is replaced by $\pi_N\varphi(\omega)$
with uniform bounds that do not depend on $N$).

\hfill$\Box$

Next we present a modified version of Lemma \ref{degreg}
that will be useful to prove Theorem \ref{main}
for $k=2$ (i.e. $m=0$ following the notation introduced in this section).
\begin{remark}\label{hopeult}
Indeed the main difference between the case $m=0$ and $m>0$ is that
Proposition \ref{carpi} is not available for $m=0$.
\end{remark}
\begin{lem}\label{degregnew}
Let  $p_3(u)\in {\mathcal P}_3(u)$ be such that
$\tilde p_3(u)=u^2 \partial_x u.$
For every $R>0$ there exists $C=C(R)>0$ such that
\begin{equation}\label{first_bis}
\mu_{1}\{u\in H^{1/2 -\epsilon} | 
|f_{N}^{p_3}(u)|>\lambda, \|\pi_N u\|_{L^2}<R, |E_{1/2}(\pi_N u)-\alpha_N|<R\}
\end{equation}$$\leq C e^{-\frac{\lambda^2}{C \alpha_N^2}},
 \hbox{ } \forall N\in \N, \lambda>0.
$$
\end{lem}
{\bf Proof.} 
First notice that due to \eqref{Runmezzo}
we have the following estimate:
$$\|\pi_N u\|_{\dot H^{1/2}}^2\leq |R_{1/2} (\pi_N u)|
+ |E_{1/2} (\pi_N u)|$$
$$
\leq C R^{3/2}\|\pi_N u\|_{H^\frac 12}^{3/2}
+ \alpha_N + R$$
$$\forall u \in \{u\in H^{1/2 -\epsilon}|\|\pi_N u\|_{L^2}<R, |E_{1/2}({\pi_N} u)-\alpha_N|<R\}.
$$
The estimate above implies
\begin{equation}\label{newtwo}
\|\pi_N u\|_{H^\frac 12}^2 \leq C(\alpha_N + 1)
\end{equation}
$$\forall u \in \{u\in H^{1/2 -\epsilon}|\|\pi_N u\|_{L^2}<R, |E_{1/2}({\pi_N}u)-\alpha_N|<R\}
$$
where $C>0$ is a suitable constant.
By combining Lemma \ref{cubichard} with \eqref{newtwo} we get
$$\Big | \int p_3(\pi_N(u)) dx\Big|\leq C \|\pi_N u\|_{H^{1/2}}^2\|\pi_N u\|_{W^{\epsilon,p}}$$
$$\leq C (\alpha_N +1) \|\pi_N u\|_{W^{\epsilon,p}}.$$
The proof can be concluded as in Lemma 6.3.

\hfill$\Box$

In the sequel the sets 
$A_{M,N}^{p, \lambda}$ and $B_{M,N}^\lambda$ 
are the ones introduced in \eqref{Amn}
and \eqref{Bmn}.
\begin{lem}\label{colsp} Let $m\geq 0$ be an integer
and $p_3(u)\in {\mathcal P}_3(u)$
such that
$\tilde p_3(u)=u\partial_x^m u \partial_x^{m+1}u.$
Then 
\begin{multline}\label{second}
\mu_{m+1}\{u\in H^{m+1/2 -\epsilon} ||f_{N}^{p_3}(u)|>\lambda,
|h_{N}(u)|<R\} 
\leq 
d\mu_{m+1} (B_{M,N}^S) 
\\
+ 
\mu_{m+1}\{u\in H^{m+1/2 -\epsilon}| 
|f_{M}^{p_3}(u)|>\frac \lambda 2
, |h_{M}(u)|<R+S\}
+
d\mu_{m+1}(A_{M,N}^{p_3, \frac \lambda 2}), \\
\forall\,\, M,N,\lambda, R, S.
\end{multline}
\end{lem}
{\bf Proof.} 
We have the following elementary estimates: 
$$
\mu_{m+1}\{u\in H^{m+1/2 -\epsilon} | |f^{p_3}_{N}(u)|>\lambda,  
|h_{N}(u)|<R\}
$$$$\leq \mu_{m+1}\{u\in H^{m+1/2 -\epsilon} | 
|f^{p_3}_{M}(u)|>\frac \lambda 2,
|h_{N}(u)|<R\}
$$$$+ \mu_{m+1}(A_{M,N}^{p_3, \frac \lambda 2})
$$
$$\leq
\mu_{m+1}\{u\in H^{m+ 1/2 -\epsilon}\setminus 
B_{M,N}^S |  |f_{M}^{p_3}(u)|>
\frac \lambda 2,
|h_{N}(u)|<R\}
$$$$+\mu_{m+1}( B_{M,N}^S)+\mu_{m+1}(A_{M,N}^{p_3, \frac \lambda 2}).
$$
On the other hand
$$\{u\in H^{m+ 1/2 -\epsilon} \setminus 
B_{M,N}^S | |h_{N}(u)|<R\}
\subset \{u \in H^{m+ 1/2 -\epsilon} |  |h_{M}(u)|<R+S\}$$
and hence we get \eqref{second}.

\hfill$\Box$

Next we propose a modified version of Lemma~ \ref{colsp} that will be useful to 
prove Theorem \ref{main} for $k=2$  (i.e. $m=0$ following the notation introduced in this section).
See Remark~\ref{hopeult} to understand the difference 
between the case  $m=0$ and $m>0$. 
\begin{lem}\label{new3} 
Let $p_3(u)\in {\mathcal P}_3(u)$
such that
$\tilde p_3(u)=u^2 \partial_xu.$
Then 
\begin{equation*}
\mu_{1}\{u\in H^{1/2 -\epsilon} ||f_{N}^{p_3}(u)|>\lambda,
\|\pi_N u\|_{L^2}<R,
|E_{1/2}(\pi_N u)-\alpha_N|<R\} 
\end{equation*}
$$\leq 
\mu_{1} \{u\in H^{1/2 -\epsilon}| |E_{1/2}(\pi_N u) -\alpha_N -E_{1/2}(\pi_M u)
+\alpha_M|>S \} 
$$$$+ 
\mu_{1}\{u\in H^{\frac 12 -\epsilon}| 
|f_{M}^{p_3}(u)|>\frac \lambda 2
, \|\pi_M u\|_{L^2}<R+S,|E_{1/2}(\pi_M u)-\alpha_M|<R+S\}$$$$+
\mu_{1}(A_{M,N}^{p_3, \frac \lambda 2}) 
+\mu_1\{u\in H^{1/2 -\epsilon}| \|\pi_M u -\pi_N u\|_{L^2}>S\},
\hbox{ }  
\forall M,N,\lambda, R, S.
$$ 
\end{lem}
{\bf Proof.} 
It is similar to the proof of Lemma \ref{colsp}.

\hfill$\Box$

{\bf Proof of Proposition \ref{proprImain45}} We have to prove \eqref{incorporame}. 
\\
\\
{\bf Claim}
\\
\\{\em It is sufficient to prove \eqref{incorporame}
with $f_{N}$ replaced by $f_{N}^{p_3}$
where } 
\begin{equation}\label{reves}
\tilde p_3(u)=u\partial_x^m u \partial_x^{m+1}u.
\end{equation}
To prove the claim first notice that due to the factor
$\prod_{j=0}^{2m} \psi(E_{j/2}(\pi_N u))$ 
in \eqref{incorporame} and due to
\eqref{carpiclaim}, we deduce that the
$L^q(d\mu_{m+1})$ norm in \eqref{incorporame} can be computed
on a sub-region $\Omega_N\subset H^{m+1/2 -\epsilon}$ such that
\begin{equation}\label{Omegan}
\Omega_N=\{u\in H^{m+ 1/2 -\epsilon}|\|\pi_N u\|_{H^m}\}<C\} 
\end{equation}
with $C$ that does not depend on $N$.  Next we prove the claim.
\\
\\
{\em First case: $m\geq 2$}
\\
\\
By combining Lemma \ref{rottcazztz}
and Lemma \ref{quarticstro}, and by looking at the structure of $E_{m+1}$
in \eqref{even} we deduce that
$$\sup_N \Big
\|R_{m+1}(\pi_N u) - \sum_{\substack{p_3(u) 
\in {\mathcal P}_3(u) s.t. \\\tilde p_3(u)=
u \partial_x^{m}u \partial_x^{m+1}u}} c_{2(m+1)}(p) \int p_3(\pi_N u) dx
\Big \|_{L^\infty(\Omega_N)}<\infty.$$
This implies the claim for $m\geq 2$.
\\
\\
{\em Second case: $m=0$}
\\
\\
In the case $m=0$ we have 
$$f_{N}(u)= \frac 34  \int 
(\pi_N u)^2 H(\pi_N \partial_x u) dx 
+ \frac 18 \int (\pi_N u)^4 dx$$
and hence
we get the claim due to 
the positivity of the last term.
\\
\\
{\em Third case: $m=1$}
\\
\\
By looking at the structure of $E_2$ and by using Lemma \ref{quarticstro}
we deduce that the claim follows provided that we prove
\begin{equation}
\sup_N \|e^{q|f^{p}(\pi_N u)|}\|_{L^1(\Omega_N, d\mu_2)}<\infty
\end{equation}
where $\Omega_N$ is defined in \eqref{Omegan}
and $p(u)\in {\mathcal P}_3(u)$ is such that
$\tilde p(u)=(\partial_x u)^3$. For simplicity we treat the case 
$p(u)= (\partial_x u)^3$
(the general case can be treated in a similar way).
We have the following estimate
$$|f^{p}(\pi_N u)|\leq C\|\pi_N u\|_{H^1}^2
\|\partial_x \pi_N u\|_{L^\infty}
\leq C \|\pi_N u\|_{W^{1+\epsilon,p}}$$
provided that $u\in \Omega_N$ and $\epsilon p>1$.
Hence we get
$$\{u\in \Omega_N| |f^{p}(\pi_N u)|>\lambda\}
\subset \{u\in \Omega_N| \|\pi_N (u)\|_{W^{1+\epsilon,p}}>\lambda C^{-1}\}$$
which in turn implies
$$\mu_2\{u\in \Omega_N| |f^{p}(\pi_N u)|
>\lambda\}\leq C e^{-\frac{\lambda^2}C}
$$
where we have used an adapted version of Proposition \ref{sobo2}
(i.e. we choose $\varphi(\omega)$ as in \eqref{randomized} with $k=4$ and 
the norm $W^{\epsilon, p}$ is replaced by $W^{1+\epsilon, p}$).
As a consequence we get
$$\int_{\Omega_N} e^{q|f^{p}(\pi_N u)|}d\mu_2
\leq \int_0^\infty e^{q\lambda} d\mu_2
\{u\in \Omega_N| |f^{p}(\pi_N u)|>\lambda\}
d\lambda<\infty.$$\\
\\
Next we shall prove \eqref{incorporame}
where $f_{N}$ is replaced by $f_{N}^{p_3}$
with $p_3(u)$ that satisfy \eqref{reves}, and it will complete the proof.\\
\\
We split the proof in two subcases.
\\
\\
{\em First case: $m>0$}
\\
\\
Since $\psi$ is compactly supported 
there exists $R>0$
such that
$$0\leq \psi(E_{j/2}(\pi_N u))\leq 
\chi_{\{u\in H^{m+1/2 -\epsilon} ||E_{j/2}(\pi_N u)|<R\}} 
$$$$\hbox{ a.e. (w.r.t. $\mu_{m+1}$) } u, \hbox{ } \forall j=0,...,2m
$$
and also
$$0\leq \psi(g_{N}(u))\leq 
\chi_{\{u\in H^{m+1/2 -\epsilon} | |g_{N}(u)|<R\}},
\hbox{ a.e. (w.r.t. $\mu_{m+1}$) } u
$$
where $\chi_A$ denotes in general the characteristic function of $A$.
In particular
$$0\leq \prod_{j=0}^{2m} \psi(E_{j/2}(u))
\psi(g_{N}(u))
$$$$\leq \chi_{\bigcap_{j=0}^{2m} 
\{u\in H^{m+1/2 -\epsilon}||E_{j/2}(\pi_N u)|< R\} 
\cap \{u\in H^{m+ 1/2 -\epsilon}||g_{N}(u)|<R\}}
$$
which due to Proposition \ref{carpi} implies
$$0\leq \prod_{j=0}^{2m} \psi(E_{j/2}(\pi_N u))
\psi(g_{N}(u))
\leq  \chi_{\{u\in H^{m+ 1/2 -\epsilon}|
|h_{N}(u)|<R\}}$$
for a suitable $R>0$ that can be different from the previous one.
Hence 
it is sufficient to prove
\begin{equation*}
\sup_N \int_{\{u\in H^{m+ 1/2 -\epsilon}||h_{N}(u)|<R\}} 
e^{q |f_{N}^{p_3}(u)|} d\mu_{m+1}<\infty, \hbox{ } \forall N\in \N
\end{equation*}
where  $\tilde p_3(u)$ satisfies \eqref{reves}. 
The estimate above is equivalent to:
\begin{equation*}
\sup_N \int_0^\infty \mu_{m+1}
\{u\in H^{m+ 1/2 -\epsilon}| |f_{N}^{p_3}(u)|>\lambda,
|h_{N}(u)|<R\}
e^{q\lambda} d\lambda<\infty.
\end{equation*}
In turn it follows by the following ones:
\begin{equation}\label{final1}\sup_N \int_{\sqrt N}^\infty  
\mu_{m+1}\{u\in H^{m+ 1/2 -\epsilon}| |f_{N}^{p_3}(u)|>\lambda, 
|h_{N}(u)|<R\} e^{q\lambda} 
d\lambda<\infty\end{equation}
and
\begin{equation}\label{final2}
\sup_N \int_0^{\sqrt N} \mu_{m+1}
\{u\in H^{m+ 1/2 -\epsilon}| |f_{N}^{p_3}(u)|>\lambda, 
|h_{N}(u)|<R\}e^{q\lambda} d\lambda<\infty. 
\end{equation}
By \eqref{first} we get:
$$\sup_N \int_{\sqrt N}^\infty 
\mu_{m+1}\{u\in H^{m+ 1/2 -\epsilon}| |f_{N}^{p_3}
(u)|>\lambda,
|h_{N}(u)|<R\}
e^{q\lambda} d\lambda
$$$$\leq \sup_N \int_{\sqrt N}^\infty C 
e^{-\frac{\lambda^2}{C\alpha_N^2}} e^{q\lambda} d\lambda
=C \sup_N  e^{\frac{Cq^2\alpha_N^2}{4}}\int_{\sqrt N}^\infty 
e^{-(\frac \lambda{\sqrt C \alpha_N}-\frac q2\sqrt 
C \alpha_ N)^2} d\lambda$$
$$=C \sup_N  e^{\frac{Cq^2\alpha_N^2}{4}} 
\int^\infty_{\frac{\sqrt N}{\sqrt C \alpha_N}-\frac q2 
\sqrt C \alpha_N} e^{-\lambda^2} 
\sqrt C \alpha_N d\lambda$$
which due to the bound
$e^{-r^2}<e^{-r}$ for every $r>1$ can be estimated
by
$$...\leq C^\frac 32 \sup_N  (\alpha_N e^{\frac{Cq^2\alpha_N^2}{4}} 
e^{-\frac{\sqrt N}{\sqrt C \alpha_N}
+\frac q2 \sqrt C \alpha_N})$$
which implies \eqref{final1}.\\
In order to prove \eqref{final2}
we use \eqref{second}
where we fix
$M=[\lambda]^2$
(here $[\lambda]$ is the integer part of $\lambda$)
and $S>0$ will be chosen later 
in a suitable way.
By recalling also \eqref{third}, \eqref{fourth} and \eqref{first}
we get:
$$\int_0^{\sqrt N} \mu_{m+1}\{u\in H^{m+ 1/2 -\epsilon}| |f_{N}^{p_3}(u)|>\lambda,
|h_{N}(u)|<R\}
e^{q\lambda} d\lambda
$$$$\leq
\int_0^{\sqrt N}  C e^{-\frac{\lambda^2}
{C\alpha_{[\lambda]^2}^2}+q\lambda}
+ 
C(e^{-\frac{[\lambda]S}{C}+q\lambda}+ 
e^{-\frac 1C (\frac \lambda 2 [\lambda])^\frac 23+q\lambda})d\lambda
$$
where
$C>0$ are uniform constant that can change at each step. Notice that if we choose $S$ large 
enough compared with
$q$ then 
$$\sup_N \int_0^{\sqrt N}  C 
e^{-\frac{\lambda^2}{C\alpha_{[\lambda]^2}^2}+q\lambda}
+ C(e^{-\frac{[\lambda]S}{C}+q\lambda}
+ e^{-\frac 1C (\frac \lambda 2 [\lambda])^\frac 23+q\lambda})d\lambda
$$
$$\leq 
\int_0^{\infty}  C e^{-\frac{\lambda^2}{C\alpha_{[\lambda]^2}^2}+q\lambda}
+ C(e^{-\frac{[\lambda]S}{C}+q\lambda}+ 
e^{-\frac 1C (\frac \lambda 2 [\lambda])^\frac 23+q\lambda})d\lambda
<\infty$$
which implies \eqref{final2}.
\\
\\
{\em Second case: $m=0$}
\\
\\
The main difference with the case $m>0$ is that we cannot apply Proposition~\ref{carpi}.
Due to the cut--off function $\psi$ it is sufficient to prove
\begin{equation*}
\sup_N \int_{\{u\in H^{1/2 -\epsilon}|\|\pi_N u\|_{L^2}<R,|E_{1/2}(\pi_N u)-\alpha_N|<R\}} 
e^{q |f_{N}^{p_3}(u)|} d\mu_{1}<\infty, \hbox{ } \forall N\in \N
\end{equation*}
where  $\tilde p_3(u)$ satisfies $\tilde p_3(u)=u^2\partial_x u$
and $R>0$ is a suitable constant. 
The estimate above is equivalent to:
\begin{equation*}
\sup_N \int_0^\infty \mu_{1}
({\mathcal C}_{\lambda,N})
e^{q\lambda} d\lambda<\infty
\end{equation*}
where
$${\mathcal C}_{\lambda,N}=\{u\in H^{1/2 -\epsilon}| |f_{N}^{p_3}(u)|>\lambda,
\|\pi_N u\|_{L^2}<R,|E_{1/2}(\pi_N u)-\alpha_N|<R\}$$
and in turn it follows by the following ones:
\begin{equation}\label{final1new}\sup_N \int_{N^\beta }^\infty  
\mu_{1} ({\mathcal C}_{\lambda,N})  e^{q\lambda} 
d\lambda<\infty\end{equation}
and
\begin{equation}\label{final2new}
\sup_N \int_0^{N^\beta } \mu_{1} ({\mathcal C}_{\lambda,N})
e^{q\lambda} d\lambda<\infty 
\end{equation}
where $\beta>0$ will be fixed later.
By Lemma \ref{degregnew} we get:
$$\sup_N \int_{N^\beta }^\infty 
\mu_{1} ({\mathcal C}_{\lambda,N})
e^{q\lambda} d\lambda
$$$$\leq \sup_N \int_{N^\beta }^\infty C 
e^{-\frac{\lambda^2}{C\alpha_N^2}} e^{q\lambda} d\lambda
=C \sup_N  e^{\frac{Cq^2\alpha_N^2}{4}}\int_{N^\beta }^\infty 
e^{-(\frac \lambda{\sqrt C \alpha_N}-\frac q2\sqrt 
C \alpha_ N)^2} d\lambda$$
$$=C \sup_N  e^{\frac{Cq^2\alpha_N^2}{4}} 
\int^\infty_{\frac{ N^\beta}{\sqrt C \alpha_N}-\frac q2 
\sqrt C \alpha_N} e^{-\lambda^2} 
\sqrt C \alpha_N d\lambda$$
which due to the bound
$e^{-r^2}<e^{-r}$ for every $r>1$ can be estimated
by
$$...\leq C^\frac 32 \sup_N  (\alpha_N e^{\frac{Cq^2\alpha_N^2}{4}} 
e^{-\frac{N^\beta }{\sqrt C \alpha_N}
+\frac q2 \sqrt C \alpha_N})$$
which implies \eqref{final1new} for every $\beta>0$.\\
In order to prove \eqref{final2new}
we use Lemma~\ref{new3} 
where we fix
$M=[\lambda]^{1/\beta}$
(here $[\lambda]$ is the integer part of $\lambda$)
and $S=1$.
By recalling also 
\eqref{third}, 
\eqref{firstnew}, \eqref{firstnew2} and \eqref{first_bis} we get
\begin{multline*}
\sup_N \int_0^{N^\beta} \mu_1 ({\mathcal C}_{\lambda, N})
e^{q\lambda} d\lambda
\leq \sup_N 
\int_0^{N^\beta}  \Big[C e^{-\frac{\lambda^2}
{C'\alpha_{[\lambda]^{1/\beta}}^2}+q\lambda}
+ 
\\
C(e^{-\frac 1C (\frac \lambda 2 [\lambda]^{1/(2\beta)})^{2/3}+q\lambda}
+ e^{-\frac 1C (\lambda [\lambda]^{1/(2\beta)})+q\lambda})\Big]d\lambda
\end{multline*}
where $C>0$ denote uniform constants which can change at each step. 
Notice that if we choose $\beta<1$ then we can continue the estimate as follows:
$$...\leq \sup_N   
\int_0^{\infty}  e^{-C \lambda^{1+\epsilon_0}+q\lambda }d\lambda
<\infty$$
for suitable $C, \epsilon_0>0$.
Hence we get \eqref{final2new}.

\hfill$\Box$

\begin{prop}\label{proprImainprime}
Let $\bar f(u),  \bar h(u)$ be as in Proposition \ref{proprImain}
and $\chi_R$ as in Theorem~\ref{main}.
Then 
$$\bigcup_{R>0}{\rm supp}
(\prod_{j=0}^{2m} \chi_R(E_{ j/2}(u))
\chi_R(\bar h(u)) + R_{m+1/2}(u)) 
e^{-\bar f(u)})= {\rm supp}(\mu_{m+1}).$$
\end{prop}
{\bf Proof.} Due to \eqref{finiti}
and \eqref{lastbnotleas} we get 
\begin{equation}\label{giul}
\mu_{m+1} \{u\in H^{m+ 1/2 -\epsilon}| |\bar h(u)|=\infty\}
\end{equation}
$$=\mu_{m+1} \{u\in H^{m+ 1/2 -\epsilon}| |R_{m+1/2}(u)|=\infty\}
=0$$
and
\begin{equation}\label{gimpi}
\mu_{m+1} \{u\in H^{m+1/2 - \epsilon}| e^{-\bar f(u)}=0\}=0
\end{equation}
Moreover by \eqref{lastbnotleas} we also get
\begin{equation}\label{lucr}
\mu_{m+1} \{u\in H^{m+ 1/2 - \epsilon}
| |E_{j/2}(u)|=\infty\}=0, \hbox{ } \forall j=0,..., 2m.
\end{equation}
As a consequence of \eqref{giul} and \eqref{lucr},
and by noticing that $\chi_R(t)\rightarrow 1$
as $R\rightarrow \infty$, we deduce:
$$\lim_{R\rightarrow \infty} 
\prod_{j=0}^{2m} \chi_R(E_{j/2}(u))
\chi_R(\bar h(u) +  R_{m+1/2}(u))=1$$$$
\hbox{ a.e. (w.r.t. $\mu_{m+1}$) } u\in H^{m+\frac 12 - \epsilon}$$ 
and hence by the Egoroff Theorem we get
$$\forall \delta>0 \hbox{ } \exists \Omega_\delta
\subset H^{m+ 1/2 - \epsilon}, \tilde R>0
\hbox{ s.t. }$$
$$ \mu_{m+1}(\Omega_\delta)> 1 -\delta \hbox{ and }
\prod_{j=0}^{2m} \chi_R(E_{j/2}(u))
\chi_R(\bar h(u) + R_{m+1/2}(u))>1-\delta$$$$
\hbox{ a.e. (w.r.t. $\mu_{m+1}$) } u\in \Omega_\delta \hbox{ and } \forall R>\tilde R.$$
By combining this fact with \eqref{gimpi} we deduce 
$$
\mu_{m+1} \Big \{\bigcup_{R>0} {\rm supp} ( \prod_{j=0}^{2m} \chi_R(E_{j/2}(u))
\chi_R(\bar h(u) + R_{m+1/2}(u))e^{-\bar f(u)})\Big\}=1.$$

\hfill$\Box$

{\bf Proof of Theorem \ref{main} for $k=2(m+1)$}
\\
It follows by combining Propositions~\ref{proprImain}, \ref{proprImain45}, \ref{proprImainprime}.

\hfill$\Box$

\section{Proof of Theorem \ref{main} for $k=2m+1$, $m>0$}

In this section we briefly describe how to adapt the proof
of Theorem \ref{main} given for $k=2(m+1)$ to the case $k=2m+1$.
We do not give the details of the proofs, however we 
underline the points where they have to be modified compared to the case $k=2(m+1)$. 
The following is an adapted version of Proposition~\ref{carpi}.
\begin{prop}\label{carpiodd}
Let $k\geq 2$ be a fixed integer. Then for every $R_1, R_2>0$ 
there is $C=C(R_1,R_2)>0$ such that
\begin{equation}\label{propgianodd}
\bigcap_{j=0}^{2k-1} \{u\in H^{k-1/2}| |E_{j/2}(u)|< R_1\} 
\cap \{u\in H^{k- 1/2}| |E_k(\pi_N u )-\alpha_N|<R_2\} 
\end{equation}
$$\subset \{ u\in H^{k-\frac 12} | \|u\|_{H^{k-1/2}}<C\}\cap 
\{u\in H^{k-1/2}|\|\pi_N u\|_{\dot H^{k}}^2 - \alpha_N|<C\}, \hbox{ }
\forall N\in \N.$$
\end{prop}
{\bf Proof.}
It is similar to the proof of Proposition \ref{carpi}, hence we skip it.

\hfill$\Box$

By looking at \eqref{odd},
in analogy with our argument used to treat the conservation laws
$E_{m+1}$,
the most delicate terms to be treated 
in $E_{m+ 1/2}$ are of the type
$\int p_3(u) dx$
where
\begin{equation}\label{p3odd}p_3(u)\in {\mathcal P}_3(u) \hbox{ and }
\tilde p_3(u)= u \partial_x^m u \partial_x^mu.\end{equation}
Next we present an adapted version of Propositions \ref{wiener}
and \ref{wiener2} in the case when $p_3(u)$ satisfies
\eqref{p3odd}. 
We recall that the Sobolev spaces $H^{m- \epsilon}$
are a support for the measure $\mu_{m+1/2}$
for every $\epsilon>0$. This fact will be used in the sequel
without any further comment.
\begin{prop}\label{wienerodd} Let $m\geq 1$ be a given integer
and 
$p_3(u)\in {\mathcal P}_3(u)$ 
be 
such that $$\tilde p_3(u)= u \partial_x^m u 
\partial_x^{m}u.$$
Then for every $\alpha\in (0, \frac 12)$
there exists $C=C(\alpha)>0$ such that 
\begin{equation}\label{2pwieodd}
\|f^{p_3}(\pi_N u)- f^{p_3}(\pi_M u)\|_{L^p(d\mu_{m+1/2})}
\leq C \frac {p^{3/2}}{(\min\{M,N\})^\alpha},  
\hbox{ } \forall M,N\in \N, p\geq 2
\end{equation}
where $f^{p_3}(v)= \int p_3(v) dx$.
In particular 
\begin{equation}\label{thirdodd}
\exists C>0 \hbox{ s.t. } 
 \mu_{m+1/2}(A_{M,N}^{p_3,\lambda})\leq  
e^{-\frac1C (\lambda \min \{N,M\}^\alpha)^{2/3}}, \hbox{ }
\forall M,N\in \N, \lambda>0
\end{equation}
where 
\begin{equation}\label{Amnodd}
A_{M,N}^{p_3,\lambda}=\{u\in H^{m-\epsilon}
||f^{p_3}(\pi_N u)- f^{p_3}(\pi_M u)|>\lambda\}.
\end{equation}
\end{prop}

\begin{prop}\label{wiener2odd} Let $m\geq 1$ be a given integer. There exists $C>0$ such that
\begin{equation}\label{2pwie2odd}
\|h_N(\pi_N u)- h_M(\pi_M u)\|_{L^p(d\mu_{m+1/2})}
\leq C \frac {p}{\sqrt {\min\{M,N\}}},  
\hbox{ } \forall M,N\in \N, p\geq 2
\end{equation}
where
$h_K(v)= \|v\|_{\dot H^m}^2- \alpha_K \hbox{ for any }  K\in \N$.
In particular 
\begin{equation}\label{fourthodd}
\exists C>0 \hbox{ s.t. } 
\mu_{m+1/2}(B_{M,N}^\lambda)\leq  
e^{-\frac1C (\lambda \sqrt {\min \{N,M\}})},
\hbox{ }
\forall M,N\in \N, \lambda>0
\end{equation}
where
\begin{equation}\label{Bmnodd}
B_{M,N}^\lambda=\{u\in H^{m-\epsilon}
| |h_{N}(\pi_N u) - h_{M}(\pi_N u)|>\lambda
\}.\end{equation}
\end{prop}
{\bf Sketch of the Proof.}
The proof of Proposition \ref{wiener2odd} is identical to the proof
of Proposition \ref{wiener2}.
Concerning the proof of Proposition \ref{wienerodd}
notice that (following the proof of Proposition
\ref{wiener}) it is sufficient to prove \eqref{2pwieodd}
for $p=2$.
By using the parametrization \eqref{randomized}
(for $k=2m+1$)
we have to estimate
$$
\Big \|\sum_{(i,j,k)\in {\mathcal A}_M^N} \frac{1}{\sqrt {|i|} |j|^{m+\frac 12} \sqrt{|k|}}
\varphi_i(\omega) \varphi_j(\omega) \varphi_k(\omega) \Big\|_{L^2_\omega}^2
$$
where ${\mathcal A}_M^N$
is the set defined in \eqref{Astorto}.
By using Lemma \ref{orth} and arguing as in Lemma \ref{orthogonal}
we can estimate the quantity above by 
$$\sum_{(i, j)\in \Z\setminus\{0\}, i+j\neq 0,|i|>\frac M2
} \frac{1}{|i||i+j|^{2m+1}|j|}\leq \frac C{M^\alpha}$$
for every $\alpha\in (0, 1)$.
The last estimate can be deduced by looking at the argument in \cite{tz}
(see end of page 500).

\hfill$\Box$

Next we present a lemma allowing us to treat all the terms that
appear in the expression of $E_{m+1/2}$ except the ones with the structure
\eqref{p3odd} (see \eqref{odd}).
\begin{lem}\label{8990}
Let $m\geq 1$ be an integer and $p_3(u)\in {\mathcal P}_3(u)$
such that 
\begin{equation}\label{cubic1}
\tilde p_3 (u)= \prod_{i=1}^3\partial_x^{\alpha_i}
\hbox{ with }
\sum_{i=1}^3\alpha_i=2m \hbox{ and } 
1\leq \min_{i=1,2,3} \alpha_i\leq 
\max_{i=1,2,3} \alpha_i \leq m.\end{equation}
Then for every $\epsilon>0$, $p\in [1, \infty)$ such that
$\epsilon p>1$ there exists $C=C(\epsilon,p)>0$ such that:
$$\Big |\int p_3(u) dx\Big |\leq C \|u\|_{H^{m-1/2}}^2 
\| u\|_{W^{m-1+\epsilon, p}}.$$
\end{lem}

{\bf Sketch of the Proof.} We treat for simplicity 
the case $p_3(u)=\tilde p_3(u)$
(the general case can be treated by a similar argument).
Next we also assume
$\alpha_1\geq \alpha_2\geq \alpha_3$.
Notice that by an integration by parts argument we can always
reduce to the following two cases.
\\
\\
{\em First case: $\alpha_1=m, \alpha_2=m-1$, $\alpha_3=1$}
\\
\\
In this case combine Lemma \ref{LP}
with the Sobolev embedding $W^{\epsilon, p} \subset L^\infty$.
\\
\\
{\em Second case: $\alpha_1\leq m-1$}
\\
\\
In this case we combine
the Cauchy-Schwartz inequality with the Sobolev embedding
$W^{\epsilon, p} \subset L^\infty$ and we get
$$\Big |\int p_3(u) dx\Big |\leq \|u\|_{H^{m-1}}^2
\|\partial_{x}^{m-1}u\|_{W^{\epsilon, p}}.$$

\hfill$\Box$

\begin{lem}\label{perpartprima}
Let $m\geq 1$ be an integer and 
$p_4(u)\in {\mathcal P}_4(u)$ such that
\begin{equation}\label{quartic1}
\tilde p_4 (u)= u^2 \partial_x^{m} u
\partial_x^{m-1}u .\end{equation}
Then for every $\epsilon>0$, $p\in [1, \infty)$ such that 
$\epsilon p>1$ there exists $C=C(\epsilon,p)>0$ such that:
$$\Big |\int p_4(u) dx\Big |\leq C \|u\|_{H^{m-1/2}}(
\|u\|_{H^1}^2\|u\|_{W^{m-1+\epsilon, p}}
+\|u\|_{H^{m-1/2}} \|u\|_{H^1}^2).
$$ 
\end{lem}

{\bf Proof.} We treat for simplicity the case
$p_4(u)=u^2 \partial_x^{m} u
\partial_x^{m-1} u$. Indeed in this specific case we could get the 
estimate as a direct application Lemma \ref{LP}.
However we propose a different and more robust 
proof that can be generalized for any $p_4(u)$
as in the assumptions.
We start by the following inequality:
\begin{equation}\label{raffin}
\Big |\int (\partial_x v_1 )v_2 dx\Big |\leq
\|v_1\|_{H^{1/2}} \|v_2\|_{H^{1/2}}.\end{equation}
By using the above estimate, where we choose 
$v_1=\partial_x^{m-1} u$ and $v_2=u^2\partial_x^{m-1}u$, 
in conjunction with the following one:
\begin{equation}\label{algquas}
\|w_1 w_2\|_{H^{1/2}}\leq C (\|w_1\|_{H^{1/2}}\|w_2\|_{L^\infty}
+ \|w_1\|_{L^\infty}\|w_2\|_{H^{1/2}})
\end{equation} we get
\begin{equation*}
\Big |\int p_4(u) dx\Big | \leq C
\|\partial_x^{m-1}u\|_{H^{1/2}} \|u^2 \partial_x^{m-1} u\|_{H^{1/2}}
\end{equation*}
$$\leq C \|u\|_{H^{m-1/2}} (\|u \partial_x^{m-1}\|_{H^{1/2}} \|u\|_{L^\infty}
+ \|u\partial_x^{m-1}u\|_{L^\infty}\|u\|_{H^{1/2}} ).$$
By using again \eqref{algquas} in conjunction with the Sobolev embedding
$W^{\epsilon, p}\subset L^\infty$ 
we get the result.

\hfill$\Box$

\begin{lem}\label{87390}
Let $m\geq 3$ be an integer and $p_j(u)\in {\mathcal P}_j(u)$
such that 
\begin{equation}\label{quintic1}
\tilde p_j (u)= \prod_{i=1}^j\partial_x^{\alpha_i} u
\hbox{ with }
\sum_{i=1}^j \alpha_i\leq 2m-1 \hbox{ and } 
\max_{i=1,...,j} \alpha_i\leq m-1.\end{equation}
Then there exists $C>0$ such that:
$$\Big |\int p_j(u) dx\Big |\leq C \|u\|_{H^{m-1}}^j.$$
\end{lem}

{\bf Sketch of the Proof.} We suppose $p_j(u)= 
\prod_{i=1}^j\partial_x^{\alpha_i} u$
with $\alpha_1\geq...\geq \alpha_j$
(the general case works with a similar argument).
By integration by parts we can reduce to two cases.
\\
\\
{\em First case: $\alpha_1=m-1, \alpha_2=m-1$}
\\
\\
In this case by assumption we get
$\alpha_i\leq 1$ for every $i=3,...,j$.
Hence by using the Cauchy-Schwartz inequality 
and the Sobolev embedding $H^1\subset L^\infty$
we get
$$\Big |\int p_{j}(u) dx\Big | \leq \|\partial_x^{m-1} u\|_{L^2}^2
\prod_{i=3}^j \|u\|_{H^2}.$$
\\
\\
{\em Second case: $\alpha_i\leq m-2 \hbox{ } \forall i=1,...,j$}
\\
\\
By using the Sobolev embedding
$H^1\subset L^\infty$ we get
$$\Big |\int p_j(u) dx\Big |\leq 
\prod_{i=1}^j\|\partial_x^{\alpha_i}u\|_{H^1}$$
and hence we conclude.

\hfill$\Box$

Next we give an adapted version of Lemma \ref{degreg}.
Recall that the functions $h_N(u)$ and $f^{p_3}(u)$ are the ones introduced 
in Propositions \ref{wienerodd}
and \ref{wiener2odd}.
\begin{lem}\label{degregodd}
Let $m\geq 1$ be an integer and $p_3(u)\in {\mathcal P}_3(u)$ be such that
$$\tilde p_3(u)=u\partial_x^mu \partial_x^{m}u.$$
For every $R>0$ there exists $C=C(R)>0$ such that
\begin{equation}\label{firstodd}
\mu_{m+1/2}\{u\in H^{m-\epsilon} | 
|f^{p_3}(\pi_N u)|>\lambda, |h_{N}(\pi_N u)|<R\}
\leq C e^{-\frac{\lambda^2}{C \alpha_N^2}}\end{equation}
$$\forall N\in \N, \lambda>0.
$$
\end{lem}
{\bf Sketch of the Proof.}
We have the following inequality
$$\Big |\int p_3(\pi_N u) dx\Big |\leq C \|\pi_N u\|_{W^{\epsilon, p}}\|
\pi_N u\|_{H^m}^2
\leq C (\alpha_N + R)\|\pi_N u\|_{W^{\epsilon, p}}$$
provided that $u$ belongs to the region on l.h.s. of
\eqref{firstodd}.
The proof can be concluded by using the following estimate
$$\mu_{m+1/2} \{u\in H^{m-\epsilon}| 
\|u \|_{W^{\epsilon, p}}
>\lambda \} \leq  C e^{-\frac{\lambda^2}{C}},
\hbox{ } \forall \lambda>0$$
whose proof is similar to the proof 
of Proposition \ref{sobo2}
(the unique difference is to use along the proof
the random vector 
$\varphi(\omega)= \sum_{n\neq 0}
\frac{\varphi_n(\omega)}{|n|^{m+1/2}}e^{{\bf i}nx}$
instead of $\varphi(\omega)= \sum_{n\neq 0}
\frac{\varphi_n(\omega)}{|n|^{m+1}}e^{{\bf i}nx}$).

\hfill$\Box$

The following version of 
Proposition \ref{proprImain}
can be easily proved. Hence we skip its proof.
\begin{prop}\label{proprImainodd}
Let $m\geq 1$ and $\psi\in C_c(\R)$ be given. 
Then there exist two
functions $\bar h(u), \bar f(u)$
measurable with respect to $\mu_{m+1/2}$
such that: 
$$|\bar h (u)|, |\bar f(u)|<\infty, \hbox{ a.e. } u\in H^{m-\epsilon};$$
$$\prod_{j=0}^{2m-1} \psi(E_{j/2}(\pi_N u))
\psi( E_{m}(\pi_N u) - \alpha_N )e^{-R_{m+1/2}(\pi_N u)}
$$$$ \hbox{ converges in measure to }$$
$$\prod_{j=0}^{2m-1} \psi(E_{j/2}(u))
\psi(\bar h(u) + R_m(u)) 
e^{-\bar f(u)}.$$
Moreover
$$|E_{j/2}(u)|, |R_m(u)|<\infty, \hbox{ a.e. (w.r.t. $\mu_{m+1/2}$) }
u\in H^{m-\epsilon}$$
$$\forall j=0,...,2m-1.$$
\end{prop}
The proof of Proposition \ref{proprImainprime}
can be easily adapted to give the following result.
\begin{prop}\label{proprImainprimeodd}
Let $\bar f(u),  \bar h(u)$ 
be as in Proposition \ref{proprImainodd}
and $\chi_R$ as in Theorem \ref{main}.
Then 
$$\bigcup_{R>0}{\rm supp}
(\prod_{j=0}^{2m-1} \chi_R(E_{j/2}(u))
\chi_R(\bar h(u)) + R_m(u)) 
e^{-\bar f(u)})= {\rm supp}(\mu_{m+1/2}).$$
\end{prop}
The last step we need in order to prove Theorem \ref{main}
in the case $k=2m+1$ is 
the following version of  
Proposition \ref{proprImain45}.
\begin{prop}\label{proprImain45odd}
Let $m\geq 1$ and $\psi\in C_c(\R)$ be given. 
For every $q\in [1,\infty)$ we have
\begin{equation}\label{incorporame*odd}
\sup_N \Big \|\prod_{j=0}^{2m-1} \psi(E_{j/2}(\pi_N u))
\psi(E_{m}(\pi_N u) -\alpha_N )e^{-R_{m+1/2}(\pi_N u)}
\Big \|_{L^q(d \mu_{m+1/2})}<\infty.
\end{equation}
\end{prop}

{\bf Sketch of the Proof.}
\\
\\
{\em First case: $m\geq 2$}
\\
\\
Arguing as in the proof of Proposition
\ref{proprImain45}  (as in the case $m>0$) and
by using Lemma \ref{degregodd}, Proposition
\ref{wienerodd} and \ref{wiener2odd} we can prove
\eqref{incorporame*odd} provided
that $R_{m+1/2}$ is replaced by $f^{p_3}$
with $p_3(u)$ that satisfy \eqref{p3odd}.
\\
\\
{\em Second case: $m=1$}
\\
\\
Arguing as in the proof of Proposition
\ref{proprImain45}  (as in the case $m=0$) 
and by using an adapted version of Lemma 
\ref{degregodd} (in the same spirit as Lemma \ref{new3})
we can prove
\eqref{incorporame*odd} provided
that $R_{1+1/2}$ is replaced by $f^{p_3}$
with $p_3(u)$ that satisfy \eqref{p3odd} with $m=1$.
\\
\\
Hence the proof of \eqref{incorporame*odd} follows 
provided that we prove
the following claim (see the analogous claim stated
along the proof of Proposition \ref{proprImain45}).
\\
\\
{\bf Claim}
\\
\\
{\em  It is sufficient to prove \eqref{incorporame*odd}
with $R_{ m+1/2}(\pi_N u)$ replaced by $f^{p_3}(\pi_N u)$
where} 
\begin{equation}\label{revesodd}
\tilde p_3(u)=u\partial_x^m u \partial_x^{m}u.
\end{equation}
Due to Proposition \ref{carpiodd}
and due to the cut-off function $\psi$ we deduce that 
the $L^q$ norms (that appear in \eqref{incorporame*odd}) are actually
computed 
in the region
$\Omega_N$ given by the condition  
\begin{equation}\label{omegaN}
\Omega_N=\{u\in H^{m-\epsilon}|
\|\pi_N u\|_{H^{m-1/2}}<C\}
\end{equation} where $C>0$ is independent on $N$.\\
Next we prove the claim.
\\
\\
{\em First case: $m\geq 3$}
\\
\\
By looking at \eqref{odd} it is sufficient to prove that
\begin{equation}\label{finalwie}
\sup_N \|e^{q|f^{p_j}(\pi_Nu)|}\|_{L^1(\Omega_N, d\mu_{m+1/2})}<\infty
\end{equation}
where $p_j(u)$ satisfy \eqref{cubic1}, \eqref{quartic1}
and \eqref{quintic1}.
Notice that if $p_j(u)$ satisfies \eqref{quintic1}
then in the region $\Omega_N$ (see \eqref{omegaN}) we get
$\sup_N \|f^{p_j}(\pi_N u)\|_{L^\infty(\Omega_N)}<\infty$
(where we have used Lemma \ref{87390})
and hence we deduce \eqref{finalwie}.
Next we treat the case when $p_3(u)$ satisfies \eqref{cubic1}.
In this case by Lemma \ref{8990} we get
$|f^{p_3}(\pi_N u)| \leq C \|\pi_N u\|_{W^{m-1+\epsilon, p}}$,
provided that $u \in \Omega_N$.
In particular
$$\{u\in \Omega_N||f^{p_3}(\pi_N u)|>\lambda \}
\subset \{ u\in \Omega_N | \|\pi_N u\|_{W^{m-1+\epsilon,p}}>\lambda C^{-1} \}$$
and hence (by using a suitable version
of Proposition \ref{sobo2})
$$\mu_{m+1/2} \{u\in \Omega_N||f^{p_3}(\pi_N u)|
>\lambda \})\leq C e^{-\frac{\lambda^2}C}$$
for a suitable $C>0$.
As a consequence we get 
$$\sup_N \int_{\Omega_N} e^{q|f^{p_3}(\pi_N u)|} d \mu_{m+1/2}
$$$$=\sup_N \int_0^\infty e^{q\lambda}  d\mu_{m+1/2} \{u
\in \Omega_N||f^{p_3}(\pi_N u)|>\lambda \}
d\lambda<\infty.$$
With a similar argument we can prove
$\sup_N \|e^{q |f^p(\pi_N u)|}\|_{L^1(\Omega_N, d\mu_{m+1/2})}<\infty$
with $p(u)$ as in  \eqref{quartic1}.
\\
\\
{\em Second case: $m=1$}
\\
\\
Looking at the structure of $E_{3/2}$
(see the introduction) we have to show that
$\sup_N \|e^{q|f^p(\pi_N u)|}\|_{L^1(\Omega_N, d\mu_{3/2})}<\infty$
where $p(u)=u^3 H \partial_x u$, $p(u)=u^2H(u\partial_x u)$, $p(u)=u^5$.\\
Notice that by the Sobolev embedding $H^{1/2} \subset L^5$
we get $$|\int (\pi_N u)^5 dx|\leq C\|\pi_N u\|_{H^{1/2}}^5<C,
\hbox{ } \forall u\in \Omega_N$$
(see \eqref{omegaN} for $m=1$) and hence we get the desired bound
when $p(u)=u^5$.\\
Next we treat the term $p(u)=u^3H\partial_xu$
(the term $p(u)=u^2 H(u\partial_xu)$ can be treated in a similar way).
By using \eqref{raffin} (used along the proof of Lemma \ref{perpartprima})
in conjunction with the estimate
$$\|v_1v_2\|_{H^{1/2}}\leq C (\|v_1\|_{H^{1/2}}\|v_2\|_{L^\infty}
+ \|v_2\|_{H^{1/2}}\|v_1\|_{L^\infty})$$ we get
$$\Big |\int (\pi_N u)^3 (H\partial_x \pi_N u)dx\Big |\leq
C\|\pi_N u\|_{H^{1/2}}^2\|\pi_N u\|_{L^\infty}^2
\leq C \|\pi_N u\|_{L^\infty}^2,
\hbox{ } \forall u\in \Omega_N$$ and
hence by Sobolev embedding $H^\frac 23\subset L^\infty$
we get 
$$...\leq C \|\pi_N u\|_{H^{2/3}}^2\leq C
\|\pi_N u\|_{H^{1/2}}\|\pi_N u\|_{H^{5/6}},
\hbox{ } \forall u\in \Omega_N.
$$
Then we deduce
$$\{u\in \Omega_N||f^p (\pi_N u)|>\lambda\}
\subset \{u\in \Omega_N| \|\pi_N u\|_{H^{5/6}}>\lambda C^{-1}\}$$
and hence (by using a suitable version
of Proposition \ref{sobo2})
$$\mu_{3/2}\{u\in \Omega_N||f^p (\pi_N u)|>\lambda\}\leq C e^{-\frac{\lambda^2}{C}}.$$
In particular
$$\sup_N \|e^{q|f^p(\pi_N u)|}\|_{L^1(\Omega_N, d\mu_{3/2})}\leq
\sup_N \int_0^\infty e^{q\lambda}e^{-\frac{\lambda^2}{C}} d\lambda<\infty.$$
\\
\\
{\em Third case: $m=2$}
\\
\\
The bound $\sup_N \|e^{q |f^{p}(\pi_N u)|}\|_{L^1(\Omega_N, d\mu_{5/2})}<\infty$ 
follows by Lemma \ref{8990} in conjunction with a suitable version
of Proposition \ref{sobo2} (used in the same spirit as above),
in the case 
$p(u)\in {\mathcal P}_3(u)
$ (but $p(u)$ does not satisfy \eqref{revesodd}).
If $p(u)\in {\mathcal P}_4(u)$ satisfies \eqref{quartic1}
then we can conclude by using Lemma \ref{perpartprima}
in conjunction with a suitable version of Proposition 
\ref{sobo2}.\\
Next we treat the case
$$p_{4}(u)\in {\mathcal P}_4(u)
\hbox{ such that } \tilde p_4(u)=u(\partial_x u)^3.$$
By combining the H\"older inequality with the Sobolev embedding $W^{\epsilon,p}
\subset L^\infty$ (provided that $\epsilon p>1$) we get
$$\Big |\int p_4(\pi_N u) dx
\Big |\leq C \|\pi_N u\|_{H^1}^2
\|\partial_x \pi_N u\|_{L^\infty}\|\pi_N u\|_{L^\infty}$$$$
\leq C \|\pi_N u\|_{H^1}^3
\|\partial_x \pi_N u\|_{W^{\epsilon, p}}\leq C 
\|\pi_N u\|_{W^{1+\epsilon,p}},
\hbox{ } \forall u\in
\Omega_N.$$
Hence we can conclude as in the previous cases
by using a suitable version of Proposition
\ref{sobo2}.\\
Finally notice that by using the Sobolev embedding $H^1\subset L^\infty$
for any $p\in [1, \infty)$
we get $\sup_N e^{q|f^p(\pi_N(u)|}<\infty$ 
in the cases $p(u)\in {\mathcal P}_5(u)$ and
$\tilde p(u)=(\partial_x u)^2 u^3$, 
$p(u)\in 
{\mathcal P}_6(u) $ and $\tilde p(u)=u^5 \partial_x u $,
$p(u)=u^7$.\\ 
The proof of the claim is concluded. 

\hfill$\Box$
 
\section{Computation of $\frac d{dt}E_{m+1}(\pi_N u(t,x))$}\label{pNstar}

In this section we shall use the notations introduced 
in Section \ref{conserlawsbo}. 
Our aim is to
construct for every $N\in \N$ and for every fixed $m\in \N$ 
a function $$G_{m+1, N}: 
{\mathcal T}_N \rightarrow \R$$
where
\begin{equation}\label{trigNpol}
{\mathcal T}_N=\Big \{\sum_{|j|\in (0,N]} c_j e^{{\bf i}jx}
| \bar c_j=c_{-j}\Big \}
\end{equation}
and such that
$$\frac d{dt}E_{m+1}(\pi_N u(t, x))=G_{m+1, N}(\pi_N u(t,x))$$
where $u(t,x)$ are solutions to the truncated 
Benjamin-Ono equation \eqref{truncBO}.\\
First we introduce some preliminary notations.
\\
\\
To $p(u)\in \cup_{n=2}^\infty {\mathcal P}_n(u)$ we associate  
a new object dependent on $N\in \N$ that will be denoted 
by $p^*_N(u)$.
\\
Let $p(u) \hbox{ be such that }$ $$\tilde p(u)=\prod_{i=1}^n 
\partial_{x}^{\alpha_i} u$$ for suitable
$0\leq \alpha_1\leq ... \leq \alpha_n$ and $\alpha_i\in \N$.
First we define $p_{i,N}^*(u)$ as the function 
obtained by $p(u)$ replacing
$\partial_x^{\alpha_i}(u)$ by $\partial_x^{\alpha_i}
(\pi_{>N} (u\partial_x u))$, i.e.
\begin{equation} 
p_{i,N}^*(u) = p(u)_{|\partial_x^{\alpha_i} u= 
\partial_x^{\alpha_i} (\pi_{>N} (u\partial_x u))}, 
\hbox{ } \forall i=1,..,n
\end{equation}
where 
$$\pi_{>N} (\sum c_j e^{{\bf i}jx})= \sum_{|j|>N}  
c_j e^{{\bf i}jx}.$$
We now define
$p_N^*(u)$ as follows:
$$p^*_N(u)=\sum_{i=1}^n p_{i,N}^*(u).$$
\begin{example}
In order to clarify the definition of $p^*_N(u)$
we give an example.\\
Assume $$p(u)= \partial_x^\alpha u 
H(\partial_x^\beta u(H\partial_x^\gamma u))$$
then 
$$p^*_N(u)= \partial_x^\alpha (\pi_{>N} (u\partial_x u)) 
H(\partial_x^\beta u(H\partial_x^\gamma u))
$$$$+ \partial_x^\alpha u H(\partial_x^\beta 
(\pi_{>N} (u\partial_x u))(H\partial_x^\gamma u)) +
\partial_x^\alpha u H(\partial_x^\beta u(H\partial_x^\gamma 
(\pi_{>N} (u\partial_x u)))).$$
\end{example}
\begin{remark}
Notice that if $p(u)\in {\mathcal P}_n(u)$ (i.e. $p(u)$ 
is homogeneous of order $n$
w.r.t. $u$)  then $p_N^*(u)$ is a function homogeneous 
of order $n+1$ for every $N\in \N$.
\end{remark}
We are now able to describe the function $G_{m+1,N}$ (see \eqref{G}).
\begin{prop}\label{defG}
For every fixed integer $m\geq 0$ and for every $N\in \N $ we have:
\begin{equation}\label{ddtE}\frac d{dt}E_{m+1}(\pi_N u(t))=  
\sum_{\substack{p(u)\in {\mathcal P}_3(u) s.t. \\ 
\tilde p(u)=u\partial_x^{m} u \partial_x^{m+1}u}}
c_{2(m+1)}(p) \int p_N^*(\pi_N u(t))dx
\end{equation}$$+
\sum_{\substack{p(u)\in {\mathcal P}_{j}(u) s.t. j=3,..., 2m+4\\ 
\|p(u)\|= 2m-j+4\\ |p(u)|\leq m}} c_{2(m+1)}(p) \int p_N^*(\pi_N u(t))dx
$$
where
$u(t,x)$ solves \eqref{truncBO} and $c_{2(m+1)}(p)$
are the same constants that appear in \eqref{even}
for $k=2(m+1)$.
\end{prop}

{\bf Sketch of the Proof.} We follow \cite{zh}
(Lemma IV.3.5 page 127).\\
Let $p(u)\in {\mathcal P}_h(u))$
be such that
$\tilde p(u)= \prod_{i=1}^h \partial_x^{\alpha_i} u$.
Then
by elementary calculus
$$\frac d{dt} \int p(u(t,x)) dx= \sum_{i=1}^h 
\int p(u)_{|\partial_x^{\alpha_i} u= \partial_x^{\alpha_i} \partial_t u} dx$$
where 
$u(t,x)$ is any regular time-dependent function.
Motivated by the identity above we introduce  
$$p_{t} (u)= \sum_{i=1}^h 
p(u)_{|\partial_x^{\alpha_i} u= \partial_x^{\alpha_i} \partial_t u}.$$
By looking at the structure of $E_{m+1}$
(see \eqref{even}) we get 
\begin{equation}\label{idder}\frac d{dt } E_{m+1} u(t, x)
=  2 \int \partial_x^{m+1} u \partial_x^{m+1}\partial_tu dx
\end{equation}
$$
 + 
\sum_{\substack{p(u)\in {\mathcal P}_3(u) s.t. \\ 
\tilde p(u)=u\partial_x^{m-1} u \partial_x^{m}u}}
c_{2(m+1)}(p) \int p_t (u)dx$$
$$+
\sum_{\substack{p(u)\in {\mathcal P}_{j}(u) s.t. j=3,..., 2m+4\\ 
\|p(u)\|= 2m-j+4\\ |p(u)|\leq m}} c_{2(m+1)}(p) \int p_t(u)dx 
$$
where $u(t,x)$ is any given time dependent function.\\
Next notice that if $u(t,x)$ solves 
\eqref{truncBO} then (due to the properties $\pi_N^2=\pi_N$ and $\pi_N + \pi_{>N}=Id$)
\begin{equation}\label{truncBOmodif}
\partial_t \pi_N u + H\partial_x^2 \pi_N
u + \big((\pi_{N}u)\partial_x(\pi_{N} u)\big)=
\pi_{>N} \big((\pi_{N}u)\partial_x(\pi_{N} u)\big)
\end{equation}
and hence if we choose in \eqref{idder}
$u(t, x)=\pi_N u(t,x)$ then we can replace the derivative 
$\partial_t \pi_N u(t, x)$, that appear 
on the r.h.s. of \eqref{idder},
by the expression
$$- H\partial_x^2 \pi_N u - \big((\pi_{N}u)(\partial_x\pi_{N} u)\big)+
\pi_{>N} \big((\pi_{N}u)(\partial_x\pi_{N} u)\big).$$
Notice that if we replace $\partial_t (\pi_N u)$ by the term
$-H\partial_x^2 \pi_N u - \big((\pi_{N}u)(\partial_x\pi_{N} u)\big)$
then we get zero on the r.h.s. of \eqref{idder}
(in fact in this way we are dealing with $\pi_N u(t, x)$ as 
with an exact solution of the Benjamin-Ono equation). However the contribution that we get when we replace 
$\partial_t (\pi_N u)$ by the term  
$\pi_{>N} \big((\pi_{N}u)\partial_x(\pi_{N} u)\big)$
is not trivial
(in fact looking at \eqref{truncBOmodif}
this term  reflects how far is $\pi_N u(t,x)$ from being  
a precise solution of the Benjamin-Ono equation).\\
Hence we deduce 
\eqref{ddtE} once we notice that in the construction above 
there is no contribution coming from the quadratic part
of $E_{m+1}$. In fact 
this contribution is given by the following quantity 
$$\int \partial_x^{m+1} (\pi_N
u) \partial_x^{m+1} \pi_{>N} \big((\pi_{N}u)\partial_x(\pi_{N} u)\big)
dx$$
which is zero by orthogonality ($\pi_N u$ is localized on the $n$ modes with $|n|\leq N$ and 
$\pi_{>N} \big((\pi_{N}u)\partial_x(\pi_{N} u)\big)$ is localized in the 
complementary modes).

\hfill$\Box$
 
\section{Some algebraic identities}

The results of this section 
will be useful along the proof
of Theorem \ref{invariance}. We recall the notation
$\pi_{>N}=Id - \pi_N.$
Moreover given a function $u(x)$ we define
$$u^+=\pi_+ u\hbox{ and } u^-=\pi_- u$$
where
$\pi_+$ (resp. $\pi_-$) is the projector on 
the positive (resp. the negative) frequencies.
We recall also that
$$H\Big (\sum_{j\in \Z\setminus \{0\}} c_j e^{{\bf i} j x} 
\Big )=
-{\bf i} \sum_{j>0} c_j e^{{\bf i} j x}  + {\bf i}\sum_{j<0} 
c_j e^{{\bf i} j x} $$
and
${\mathcal T}_N$ is defined by \eqref{trigNpol}.
\begin{lem}\label{intpar1}
Let $u\in {\mathcal T}_N$ be such that $\int u dx=0$. 
Then the following identities occur:
\begin{equation}\label{mire1}
\int u (H\partial_x^m \pi_{>N} 
(u\partial_x u)) \partial_x^{m+1} u dx
\end{equation}
$$=\sum_{j=1}^m a_j [\int \pi_{>N} (\partial^j_x u^+
\partial^{m-j+1}_x u^+ )\pi
_{>N} (u^- \partial^{m+1}_xu^-)$$
$$- \pi_{>N} (\partial_x^j 
u^-\partial_x^{m-j+1} u^-)\pi_{>N} 
(u^+ \partial_x^{m+1}u^+)]$$
for suitable coefficient $a_j\in \C$;
\begin{equation}\label{mire2}
\int u (H\partial_x^m u) \partial_x^{m+1} \pi_{>N} (u \partial_x u) dx
\end{equation}$$=\sum_{j=1}^m b_j [\int \pi_{>N} 
(\partial_x^j u^+\partial_x^{m-j+1} u^+ )\pi_{>N} 
(u^- \partial_x^{m+1}u^-) $$
$$- \pi_{>N} (\partial^j_x u^-\partial^{m-j+1}_x u^-)\pi_{>N} 
(u^+ \partial^{m+1}_xu^+)]$$ for suitable coefficient $b_j\in \C$.
\end{lem}
\begin{remark}\label{usef}
Notice that 
the l.h.s. of \eqref{mire1} and \eqref{mire2} 
involve at first insight (after developing the $m$-derivative of the product)
a term that contains the product of two derivatives of order $m+1$,
which is quite dangerous (see the end of Section~\ref{kdvrie}).
The main point of the lemma above
is that on the r.h.s. of \eqref{mire1}
and \eqref{mire2} this bad term is disappeared.
\end{remark}

{\bf Proof.}
We prove \eqref{mire1}.
Due to the following identity
\begin{equation}\label{selfadj}
\int (\pi_{>N}f)g dx=\int (\pi_{>N}f)(\pi_{>N}g) dx
\end{equation}
we get:
$$\int u (H\partial_x^m \pi_{>N} (u\partial_x u)) \partial_x^{m+1} u dx
$$
$$
= \int (H\partial_x^m \pi_{>N} (u\partial_x u)) 
\pi_{>N}(u \partial_x^{m+1} u) dx.
$$
On the other hand if  $v(x), w(x)$ are 
trigonometric polynomial of degree $N$ 
we have
\begin{equation*}
\pi_{>N} (v^+ w^-)=0
\end{equation*}
and in particular
\begin{equation}\label{aselfaj}
\pi_{>N}(vw)= \pi_{>N} (v^+ w^+)+\pi_{>N} (v^- w^-).
\end{equation}
As a consequence we continue the identity above as follows
\begin{equation}\label{primi...}
...= -{\bf i} \int \partial_x^m \pi_{>N} (u^+\partial_x u^+) 
\pi_{>N} (u^- \partial_x^{m+1} u^-)dx
\end{equation}
$$
+ {\bf i}\int \partial_x^m \pi_{>N} (u^-\partial_x u^-) 
\pi_{>N} (u^+ \partial_x^{m+1} u^+)dx$$
where we have used the
definition of the Hilbert transform $H$,
$$...=
-{\bf i}\int  \pi_{>N} (u^+\partial_x^{m+1}u^+) 
\pi_{>N} (u^- \partial_x^{m+1} u^-)$$
$$+{\bf i} \int \pi_{>N} (u^-\partial_x^{m+1}u^-) 
\pi_{>N} (u^+ \partial_x^{m+1} u^+)dx$$
$${-\bf i}
\int  \pi_{>N} (\partial_x^m (u^+ \partial_x u^+) - u^
+\partial_x^{m+1}u^+) 
\pi_{>N} (u^- \partial_x^{m+1} u^-)dx$$
$$+{\bf i} \int \pi_{>N} (\partial_x^m 
(u^- \partial_x u^-) - u^-\partial_x^{m+1}u^-)
\pi_{>N} (u^+ \partial_x^{m+1} u^+)dx.$$
We can conclude by the Leibnitz rule
since the first two terms above cancel.
Concerning \eqref{mire2} notice that by using
\eqref{selfadj}
and \eqref{aselfaj} we get:
$$\int u (H\partial_x^m u) \partial_x^{m+1} \pi_{>N} (u \partial_xu) 
dx$$
$$={\bf i}\int \partial^{m+1}_x \pi_{>N} (u^+\partial_x u^+) \pi_{>N}
( u^- \partial_x^m u^-)dx
$$$$-{\bf i} \int \partial_x^{m+1} \pi_{>N} (u^-\partial_x u^-) \pi_{>N}
( u^+ \partial_x^m u^+)dx
$$
and by integration by parts
$$...=-{\bf i}\int \partial^{m}_x \pi_{>N} (u^+\partial_x u^+) 
\pi_{>N} \partial_x ( u^- \partial_x^m u^-) 
$$$$+{\bf i} \int \partial_x^{m} \pi_{>N} (u^-\partial_x u^-) 
\pi_{>N} \partial_x ( u^+ \partial_x^m u^+)dx
$$
which in turn gives
$$...=-{\bf i}\int \partial^{m}_x \pi_{>N} (u^+\partial_x u^+) 
\pi_{>N} (\partial_x u^- \partial_x^m u^-) 
$$$$+{\bf i} \int \partial_x^{m} \pi_{>N} (u^-\partial_x u^-) 
\pi_{>N} (\partial_x u^+ \partial_x^m u^+)dx
$$
$$-{\bf i}\int \partial^{m}_x \pi_{>N} (u^+\partial_x u^+) 
\pi_{>N} ( u^- \partial_x^{m+1} u^-) 
$$$$+{\bf i} \int \partial_x^{m} \pi_{>N} (u^-\partial_x u^-) 
\pi_{>N} ( u^+ \partial_x^{m+1} u^+)dx.
$$
Notice that the last two integrals above
$$-{\bf i}\int \partial^{m}_x \pi_{>N} (u^+\partial_x u^+) 
\pi_{>N} ( u^- \partial_x^{m+1} u^-) 
$$$$+{\bf i} \int \partial_x^{m} \pi_{>N} (u^-\partial_x u^-) 
\pi_{>N} ( u^+ \partial_x^{m+1} u^+)dx
$$
can be treated as in \eqref{primi...}.
Hence we have to deal with the remaining terms
in the identity above:
$$-{\bf i}\int \partial^{m}_x \pi_{>N} (u^+\partial_x u^+) 
\pi_{>N} (\partial_x u^- \partial_x^m u^-) 
$$$$+{\bf i} \int \partial_x^{m} \pi_{>N} (u^-\partial_x u^-) 
\pi_{>N} (\partial_x u^+ \partial_x^m u^+)dx.
$$
Those integrals can be easily 
handled by using the Leibnitz rule.

\hfill$\Box$

In the same spirit as in Lemma \ref{intpar1}
one can prove that if $u\in {\mathcal T}_N$ is such that $\int u dx=0$, 
then the following identities occur:
\begin{equation}\label{intpar21}
\int u (\partial_x^m \pi_{>N} (u\partial_x u)) \partial_x^{m+1} (H u) dx
\end{equation}
$$=\sum_{j=1}^m c_j 
[\int \pi_{>N} (\partial_x^j u^+\partial_x^{m-j+1} u^+ )\pi_{>N} 
(u^- \partial^{m+1}_xu^-) $$
$$- \pi_{>N} (\partial_x^j u^-\partial_x^{m-j+1} u^-)
\pi_{>N} (u^+ \partial^{m+1}_xu^+)]$$
for suitable $c_j\in \C$;
\begin{equation}\label{intpar22}
\int u \partial_x^m u (\partial_x^{m+1} H \pi_{>N} (u \partial_x u)) dx
\end{equation}
$$=\sum_{j=1}^m d_j [\int \pi_{>N} 
(\partial_x^j u^+\partial_x^{m-j+1} u^+ )
\pi_{>N} (u^- \partial_x^{m+1}u^-) $$
$$- \pi_{>N} (\partial_x^j u^-\partial_x^{m-j+1} u^-)
\pi_{>N} (u^+ \partial_x^{m+1}u^+)]$$
for suitable $d_j\in \C$;

\begin{equation}\label{intpar4567}
\int H u (\partial_x^m \pi_{>N} H (u\partial_x
u)) (\partial_x^{m+1} H u) dx
\end{equation}
$$=\sum_{j=1}^m e_j [\int \pi_{>N} 
(\partial_x^j u^+\partial_x^{m-j+1} u^+ )
\pi_{>N} (u^- \partial_x^{m+1}u^-) $$
$$- \pi_{>N} (\partial_x^j u^-\partial_x^{m-j+1} u^-)
\pi_{>N} (u^+ \partial_x^{m+1}u^+)]$$
for suitable $e_j\in \C$;

\begin{equation}\label{intpar4567prime}
\int H u (\partial_x^m H u) (\partial_x^{m+1} \pi_{>N} H (u 
\partial_x u)) dx\end{equation}
$$=\sum_{j=1}^m f_j [\int \pi_{>N} 
(\partial_x^j u^+\partial_x^{m-j+1} u^+ )
\pi_{>N} (u^- \partial_x^{m+1}u^-) $$
$$- \pi_{>N} (\partial_x^j u^-\partial_x^{m-j+1} u^-)
\pi_{>N} (u^+ \partial_x^{m+1}u^+)]$$
for suitable $f_j\in \C$;

\begin{equation}\label{intpar4567qua}
\int H u (\partial_x^m \pi_{>N} (u\partial_x
u)) (\partial_x^{m+1} u) dx
\end{equation}
$$=\sum_{j=1}^m g_j [\int \pi_{>N} 
(\partial_x^j u^+\partial_x^{m-j+1} u^+ )
\pi_{>N} (u^- \partial_x^{m+1}u^-) $$
$$- \pi_{>N} (\partial_x^j u^-\partial_x^{m-j+1} u^-)
\pi_{>N} (u^+ \partial_x^{m+1}u^+)]$$
for suitable $g_j\in \C$;

\begin{equation}\label{intpar4567quaprime}
\int H u (\partial_x^m u) (\partial_x^{m+1} \pi_{>N} (u 
\partial_x u)) dx
\end{equation}
$$=\sum_{j=1}^m h_j [\int \pi_{>N} 
(\partial_x^j u^+\partial_x^{m-j+1} u^+ )
\pi_{>N} (u^- \partial_x^{m+1}u^-) $$
$$- \pi_{>N} (\partial_x^j u^-\partial_x^{m-j+1} u^-)
\pi_{>N} (u^+ \partial_x^{m+1}u^+)]$$
for suitable $h_j\in \C$.

\begin{lem}\label{intpar3}
Let $u\in {\mathcal T}_N$ be
such that $\int u dx=0$.
Then the following identities occur:
\begin{equation}\label{intpar31}
\int u (H\partial_x^m \pi_{>N} (u\partial_x u)) 
(\partial_x^{m+1} H u) dx
\end{equation}
\begin{equation*}
+\int u (H\partial_x^m u) 
(\partial_x^{m+1} \pi_{>N} H (u \partial_x u)) dx
\end{equation*}
$$=
-\int \partial^{m}_x(\pi_{>N} (u^+\partial_x u^+))
\pi_{>N} (\partial_x u^-\partial^{m}_x u^-)dx $$$$
-
\int \partial^{m}_x(\pi_{>N} (u^-\partial_x u^-))
\pi_{>N} (\partial_x u^+\partial^{m}_x u^+)dx.
$$ 
\end{lem}
\begin{remark}
To understand the interest of Lemma \ref{intpar3},
see remark \ref{usef}.
\end{remark}
{\bf Proof.} By combining \eqref{selfadj} with 
\eqref{aselfaj} we get:
\begin{equation}\label{.}
\int u (H\partial_x^m \pi_{>N} (u\partial_xu)) 
(\partial_x^{m+1} H u) dx
\end{equation}
\begin{equation*}
+
\int u (H\partial_x^m u) (\partial_x^{m+1} H (\pi_{>N} (u 
\partial_x u))) dx
\end{equation*}
$$= \int (\pi_{>N} \partial_x^m (u^+\partial_x u^+))
\pi_{>N} (u^-\partial^{m+1}_x u^-)dx
$$
$$
+
\int (\pi_{>N} \partial_x^m (u^-\partial_x u^-))
\pi_{>N} (u^+\partial^{m+1}_x u^+)dx
$$
$$+\int \partial^{m+1}_x(\pi_{>N} (u^+\partial_x u^+))
\pi_{>N} (u^-\partial^{m}_x u^-)dx
$$
$$
+
\int (\pi_{>N} \partial_x^{m+1} (u^-\partial_x u^-))
\pi_{>N} (u^+\partial^{m}_x u^+)dx.
$$
On the other hand by integration by parts in the second term we get:
$$\int (\pi_{>N} \partial_x^m (u^+\partial_x u^+))
\pi_{>N} (u^-\partial^{m+1}_x u^-)
$$
$$+\int \partial^{m+1}_x(\pi_{>N} (u^+\partial_x u^+))
\pi_{>N} (u^-\partial^{m}_x u^-)
$$
$$=\int (\pi_{>N} \partial_x^m (u^+\partial_x u^+))
\pi_{>N} (u^-\partial^{m+1}_x u^-)
$$
$$-\int \partial^{m}_x(\pi_{>N} (u^+\partial_x u^+))
\pi_{>N} \partial_x (u^-\partial^{m}_x u^-).
$$  
By developing the derivative  
$\partial_x (u^-\partial^{m}_x u^-)= \partial_x u^-\partial^{m}_x u^-
+u^-\partial^{m+1}_x u^-$ and by replacing it in the last integral,
we get
\begin{equation}\label{...}
...=-\int \partial^{m}_x(\pi_{>N} (u^+\partial_x u^+))
\pi_{>N} (\partial_x u^-\partial^{m}_x u^-).
\end{equation}  
By using integration by parts in the second integral we get
\begin{equation}\label{....} 
\int (\pi_{>N} \partial_x^m (u^-\partial_x u^-))
\pi_{>N} (u^+\partial_x^{m+1} u^+)
\end{equation}
$$
+ 
\int \pi_{>N} (\partial_x^{m+1} (u^-\partial_x u^-))
\pi_{>N} (u^+\partial_x^{m} u^+)
$$
$$=-\int \pi_{>N} (\partial_x^{m} (u^-\partial_x u^-))
\pi_{>N} (\partial_x u^+\partial^{m}_x u^+).
$$
The proof follows by combining
\eqref{.}, \eqref{...}, \eqref{....}.

\hfill$\Box$

By a similar argument it is possible to prove that if
$u(x)$ is as in Lemma \ref{intpar3}
then the following identities occur:

\begin{equation}\label{comes}
\int H u (\partial_x^m  \pi_{>N} (u \partial_x u) )
\partial_x^{m+1} H u 
dx
\end{equation}
\begin{equation*}
+\int H u (\partial_x^m  u )\partial_x^{m+1} (\pi_{>N} H (u \partial_x u)) 
dx
\end{equation*}
$$=-
\int \partial^{m}_x(\pi_{>N} (u^+\partial_x u^+))
\pi_{>N} (\partial_x u^-\partial^{m}_x u^-)dx $$$$
-
\int \partial^{m}_x(\pi_{>N} (u^-\partial_x u^-))
\pi_{>N} (\partial_x u^+\partial^{m}_x u^+)dx.
$$ 

\begin{equation}\label{intpar23}
\int H u (\partial_x^m \pi_{>N} H (u\partial_x u)) 
\partial_x^{m+1} u dx
\end{equation}
\begin{equation*}
+\int H u (\partial_x^m H u )\partial_x^{m+1} (\pi_{>N} (u \partial_x u)) 
dx
\end{equation*}
$$=
-\int \partial^{m}_x(\pi_{>N} (u^+\partial_x u^+))
\pi_{>N} (\partial_x u^-\partial^{m}_x u^-)dx $$$$
-
\int \partial^{m}_x(\pi_{>N} (u^-\partial_x u^-))
\pi_{>N} (\partial_x u^+\partial^{m}_x u^+)dx.
$$ 

\begin{equation}\label{intpar4134}
\int u (\partial_x^m \pi_{>N} 
(u\partial_x u)) \partial_x^{m+1} u dx
\end{equation}
\begin{equation*}
+\int u (\partial_x^m u) 
\partial_x^{m+1} \pi_{>N} (u \partial_xu) 
dx\end{equation*}
$$=-
\int \partial^{m}_x(\pi_{>N} (u^+\partial_x u^+))
\pi_{>N} (\partial_x u^-\partial_x^{m} u^-)
$$$$-
\int \partial^{m}_x(\pi_{>N} (u^-\partial_x u^-))
\pi_{>N} (\partial_x u^+\partial^{m}_x u^+).
$$

\section{Some calculus inequalities}

Next we present some useful results related to the convergence
of suitable numerical series.

\begin{lem}\label{prod}
The following estimate occurs :
\begin{equation}\label{double}
\sum_{\substack{|n+m|>N\\0< |n|, |m|\leq N}} 
\frac 1{n^2}\frac 1{|m|}=O\Big (\frac{\ln N}{N}\Big )
\hbox{ as } N\rightarrow \infty.
\end{equation}
\end{lem}

{\bf Proof.}
We have the identity 
$$\sum_{\substack{|n+m|>N\\ 0< |n|, |m|\leq N}} 
\frac 1{n^2}\frac 1{|m|}
=2 \sum_{\substack{n+m>N\\ 0<n, m\leq N}} 
\frac 1{n^2}\frac 1{m}
$$
where we have used
$$\{(n,m)\in \Z \times \Z | 0<|n|, |m|\leq N, |n+m|>N \}
$$
$$=\{(n,m)\in \Z\times \Z | 0<n, m\leq N, |n+m|>N \}
$$$$\cup
\{(n,m)\in \Z \times \Z | -N \leq n, m<0, |n+m|>N \}.$$
Next we continue the identity above
$$...= 2\sum_{0<n\leq N} \frac 1{n^2} \Big (
\sum_{N-n<m\leq N}\frac 1m \Big)
\leq 2 \sum_{0<n\leq N} \frac 1{n^2} \frac n{N-n}
$$$$= \frac 2N  \sum_{0<n\leq N} \Big 
( \frac 1n + \frac 1{N-n}\Big ).$$
The proof follows since
$\sum_{0<n\leq N } \frac 1n=O(\ln N)$.

\hfill$\Box$

\begin{lem}\label{orthspa}
The following estimate occurs :
\begin{equation}\label{doublespa}
\sum_{\substack{|n+m+l|>N\\ 0< |n|, |m|, |l|\leq N}} 
\frac 1{n^2 m^2 |l|}=O\Big (\frac{\ln N}{N}\Big)
\hbox{ as } N\rightarrow \infty.
\end{equation}
\end{lem}

{\bf Proof.}
We split the sum as follows:
$$\sum_{\substack{|n+m+l|>N\\ 0< |n|, |m|, |l|\leq N}} \frac 1{n^2 m^2 |l|}
$$$$\leq \sum_{\substack{|n+l|> \frac N2\\ 
0< |n|, |m|, |l|\leq N}} 
\frac 1{n^2 m^2 |l|}
+ 
\sum_{\substack 
{|m|> \frac N2\\ 0< |n|, |m|, |l|\leq N}} 
\frac 1{n^2 m^2 |l|}
$$
$$=I_N+II_N.$$
By using Lemma \ref{prod} we get
$$I_N= O\Big( \frac{\ln N}{N}\Big).$$
Concerning $II_N$  we have
$$II_N\leq \Big (\sum_{\frac N2 < |m|\leq N} 
\frac 1{m^2}\Big)\Big( \sum_{0<|l|\leq N} 
\frac 1{|l|}\Big)\Big( \sum_{0<|n|\leq N} \frac 1{|n|^2}\Big)
\leq C \frac{\ln N}{N}.
$$

\hfill$\Box$

\section{Proof of Theorem \ref{invariance}}

Along this section we shall write 
$\varphi_N(\omega)=\pi_N(\varphi(\omega))$,
$\varphi_N^\pm (\omega)=\pi_\pm (\pi_N\varphi(\omega))$
(where $\pi_\pm$ are the projectors on the positive and negative frequencies)
and
$$\varphi(\omega)=
\sum_{n\in \Z \setminus \{0\}}
\frac{\varphi_n(\omega)}{|n|^{m+1}}e^{{\bf i}nx}.$$
Moreover for any given $p(u)\in \cup_{n=1}^\infty{\mathcal P}_n(u)$ and $N\in \N$,
$p_N^*(u)$ is defined in Section~\ref{pNstar}.\\
Notice that due to the H\"older inequality
the standard gaussian variables $\{\varphi_k(\omega)\}_k$
satisfy:
$$\forall q\in [1,\infty), k\in \N
\hbox{ } \exists C=C(k, q)>0 \hbox{ s.t. }
\sup_{j_1,...,j_k\in \Z\setminus \{0\}}
\|\varphi_{j_1}...\varphi_{j_k}\|_{L^q_\omega}\leq C.
$$
This fact will be freely used in the sequel.
\begin{lem}\label{p3sing}
Let $m\geq 2$ be an integer and $p(u)\in {\mathcal P}_3(u)$ such that
$\tilde p(u)=u\partial_x^mu \partial_x^{m+1}u.$
Then for every $q\in [1, \infty)$ we have the following
$$\lim_{N\rightarrow \infty} \Big \|\int p^*_N(\pi_N u) dx\Big \|_{L^q(d\mu_{m+1})}=0.$$  
\end{lem}

{\bf Proof.} 
By using elementary properties of the Hilbert transform 
(i.e. $H^2=-Id, \int (Hv)w dx=\int v(Hw) dx$) 
it is easy to check that if $p(u)$ is like in the assumptions then  
the quantities $\int p (u) dx$ can be always reduced 
to the following ones: 
$$\pm \int u\partial_x^mu \partial_x^{m+1}u dx, 
\pm \int (H u) (\partial_x^m H u) (\partial_x^{m+1}H u) dx,$$$$
\pm \int u(\partial_x^ mHu) (\partial_x^{m+1}Hu) dx,
\pm \int (H u)\partial_x^mu \partial_x^{m+1}u dx,
$$$$\pm \int u\partial_x^mu (\partial_x^{m+1} Hu) dx, 
\pm \int (Hu)(\partial_x^mHu) \partial_x^{m+1}u dx,$$
$$\pm \int u(\partial_x^mHu) (\partial_x^{m+1}u) dx,
\pm \int (Hu)(\partial_x^mu) (\partial_x^{m+1}Hu) dx.$$
{\em First case: $p(u)=u\partial_x^mu \partial_x^{m+1}u$}
\\
\\
In this case we can write explicitly
$$p_N^*(u)= \pi_{>N} (u \partial_x u)\partial_x^mu \partial_x^{m+1}u
$$$$+ 
u \partial_x^m (\pi_{>N} (u\partial_x u)) \partial_x^{m+1}u 
+ u\partial_x^m u \partial_x^{m+1}(\pi_{>N} (u\partial_x u))).$$
Hence we get
$$\int p_N^*(\pi_N(\varphi(\omega)))dx= 
I_N(\omega) + II_N(\omega)$$
where
\begin{equation}\label{INNN}
I_N(\omega)= \int \pi_{>N} 
(\varphi_N(\omega) \partial_x (\varphi_N(\omega)))
\partial_x^m \varphi_N(\omega)
\partial_x^{m+1}\varphi_N(\omega) dx
\end{equation}
and
\begin{equation}\label{INN9}II_N(\omega)= \int 
\varphi_N(\omega) (\partial_x^m \pi_{>N} (\varphi_N(\omega)
(\partial_x \varphi_N(\omega))) 
\partial_x^{m+1}\varphi_N(\omega)) 
\end{equation}$$+ \varphi_N(\omega)(\partial_x^m 
\varphi_N(\omega))\partial_x^{m+1}(\pi_{>N} (\varphi_N(\omega)
\partial_x(\varphi_N(\omega)))dx.$$
In order to estimate
$I_N$ notice that 
$$I_N(\omega)=\int (\pi_{>N} 
\varphi_N(\omega)\partial_x \varphi_N(\omega))
\partial_x^m \varphi_N(\omega) \partial_x^{m+1}
\varphi_N(\omega) dx
$$$$= \sum_{\substack{0<|j_1|, |j_2|, |j_3|, |j_4|\leq N\\
|j_1+j_2|>N\\j_1+j_2+j_3+j_4=0}} 
\frac{\varphi_{j_1}(\omega)}{|j_1|^{m+1}}
\frac{\varphi_{j_2}(\omega)}{|j_2|^m}
\frac{\varphi_{j_3}(\omega)}{|j_3|}\varphi_{j_4}(\omega)$$
and hence by the Minkowski inequality
$$\|I_N(\omega)\|_{L^q_\omega}\leq C 
\sum_{\substack{0<|j_1|, |j_2|, |j_3|, |j_4|\leq N\\|j_1+j_2|>N\\j_1+j_2+j_3+j_4=0}}
\frac 1{|j_1|^{m+1}|j_2|^m|j_3|}$$
$$\leq C \Big (\sum_{0<|j_3|\leq N} \frac 1{|j_3|} \Big)
\Big(\sum_{\substack{0<|j_1|, |j_2|\leq N\\|j_1+j_2|>N}}
\frac 1{|j_1|^{m+1}|j_2|^m} \Big)=O\Big(\frac{\ln^2 N}{N}\Big)$$
where we have used Lemma \ref{prod}.\\
Next we estimate $II_N(\omega)$ (see \eqref{INN9}). Due to the identity
\eqref{intpar4134}
we are reduced to estimate the following quantities:
$$II_N'(\omega)=-
\int \partial^{m}_x(\pi_{>N} (\varphi_N^+(\omega)
\partial_x \varphi_N^+(\omega)))
\pi_{>N} (\partial_x \varphi_N^-(\omega)
\partial_x^{m} \varphi_N^-(\omega))dx
$$
$$II_N''(\omega)=
-\int \partial^{m}_x(\pi_{>N} (\varphi_N^-(\omega)
\partial_x \varphi_N^-(\omega)))
\pi_{>N} (\partial_x \varphi_N^+(\omega)\partial^{m}_x 
\varphi_N^+(\omega))dx .$$
Next we estimate 
$II_N'(\omega)$ (a similar argument 
works for $II_N''(\omega)$).
By the Leibnitz formula it is sufficient to prove that:
\begin{equation}\label{jzerojenne}
\left \|\int \pi_{>N} (\partial_x^j \varphi_N^+(\omega) 
\partial^{m-j+1}_x\varphi^+_N(\omega))
\pi_{>N} (\partial_x \varphi_N^-(\omega) \partial_x^{m} 
\varphi_N^-(\omega)) dx \right \|_{L^q_\omega}=o(1)
\end{equation}
$$\hbox{ as } N\rightarrow \infty 
\hbox{ } \forall j=0,1,.., m.$$
Indeed the most delicate cases are $j=0,m$. 
All the other cases can treated in the same way.
In the case $j=0$ we are reduced to prove
$$\lim_{N\rightarrow \infty}
\left \|\int \pi_{>N} (\varphi_N^+(\omega)
\partial^{m+1}_x\varphi_N^+(\omega))
\pi_{>N} (\partial_x \varphi_N^-(\omega)\partial_x^{m} 
\varphi^-_N(\omega)) dx\right \|_{L^q_\omega}=0.$$
For that purpose, we write
$$\limsup_{N\rightarrow \infty}
\Big \|\sum_{\substack{0<|j_1|, 
|j_2|, |j_3|, |j_4|\leq N\\
j_1, j_2>0, j_3, j_4<0\\
|j_1+j_2|>N\\j_1+j_2+j_3+j_4=0}} 
\frac{\varphi_{j_1}(\omega)}{|j_1|^{m+1}} 
\varphi_{j_2}(\omega)
\frac{\varphi_{j_3}(\omega)}{|j_3|^{m}} 
\frac{\varphi_{j_4}(\omega)}{|j_4|}
\Big \|_{L^q_\omega}$$
$$\leq \limsup_{N\rightarrow \infty} 
C \sum_{\substack{0<|j_1|,|j_3|, |j_4|
\leq N\\
|j_3+j_4|>N}} \frac 1{|j_1|^{m+1} |j_3|^{m}|j_4|}
$$$$\leq C \limsup_{N\rightarrow \infty} 
\Big (\sum_{0<|j_1|\leq N} 
\frac 1{|j_1|^{m+1}} \Big) 
\sum_{\substack{0<|j_3|, |j_4|\leq N\\
|j_3+j_4|>N}} \frac 1{|j_3|^{m}|j_4|}=O
\Big (\frac{\ln N}{N}
\Big )$$
where we have used Lemma \ref{prod} at the last step.\\
To prove \eqref{jzerojenne} 
for $j=m$ we have to show
$$\lim_{N\rightarrow \infty}
\left \|\int (\pi_{>N} (\partial_x^m\varphi_N^+(\omega)
\partial_x\varphi_N^+(\omega)))
\pi_{>N} (\partial_x \varphi_N^-(\omega)\partial_x^{m} 
\varphi^-_N(\omega))\right \|_{L^q_\omega}=0.$$
Indeed arguing as above (i.e. we replace the random vector
$\varphi(\omega)$ by its random Fourier series 
and we apply the
Minkowski inequality) we are reduced to prove that
$$\lim_{N\rightarrow \infty} 
\sum_{\substack{0<|j_1|, |j_2|, |j_3|, |j_4|\leq N\\
j_1, j_2>0, j_3, j_4<0\\
|j_1+j_2|>N\\j_1+j_2+j_3+j_4=0}}
\frac 1{|j_1||j_2|^m|j_3|^m|j_4|}=0.
$$
This estimate follows
by combining the inequality 
$$\sum_{\substack{0<|j_1|, |j_2|, |j_3|, |j_4|\leq N\\
j_1, j_2>0, j_3, j_4<0\\
|j_1+j_2|>N\\j_1+j_2+j_3+j_4=0}}
\frac 1{|j_1||j_2|^m|j_3|^m|j_4|}
\leq \Big (\sum_{0<|j_3|\leq N} \frac 1{|j_3|^m}
\Big)\Big (\sum_{\substack{0<|j_1|, |j_2|\leq N\\
|j_1+j_2|>N}}
\frac 1{|j_1||j_2|^m}\Big )$$
with Lemma \ref{prod}.
\\
\\{\em Second case: $p(u)=u\partial_x^m H u \partial_x^{m+1} Hu,
H u\partial_x^m  u \partial_x^{m+1}H u, 
H u\partial_x^m H u \partial_x^{m+1} u$}
\\
\\
All those cases can be treated as the previous one
provided that we use 
\eqref{intpar31}, \eqref{comes},
\eqref{intpar23} instead of \eqref{intpar4134}
in the argument above.
\\
\\
{\em Third case: $p(u)=u\partial^m_x u (\partial_x^{m+1}Hu)$}
\\
\\
By definition we get:
$$\int p_N^*(\pi_N(\varphi(\omega)))dx
= I_N(\omega) + II_N(\omega)+III_N(\omega)$$
where
$$I_N(\omega)= \int (\pi_{>N} 
(\varphi_N(\omega)\partial_x \varphi_N(\omega)))
(\partial_x^m \varphi_N(\omega)) (\partial_x^{m+1} 
H \varphi_N(\omega)) dx,
$$
$$II_N(\omega)= \int 
\varphi_N(\omega) \partial_x^m (\pi_{>N} (\varphi_N(\omega)
\partial_x \varphi_N(\omega))) 
\partial_x^{m+1} H \varphi_N(\omega) dx
$$$$III_N(\omega)= \int 
\varphi_N(\omega)\partial_x^m 
\varphi_N(\omega) \partial_x^{m+1}H 
(\pi_{>N} (\varphi_N(\omega)
\partial_x \varphi_N(\omega))) dx.$$
The term $I_N(\omega)$ can be estimated in the same way as 
\eqref{INNN} in the first case.\\
Concerning $II_N(\omega)$ we use \eqref{intpar21}
and we get
$$II_N(\omega)
=\sum_{j=1}^m c_j [\int \pi_{>N} 
(\partial_x^j \varphi_N^+(\omega)
\partial_x^{m-j+1} \varphi^+_N(\omega))\pi_{>N} 
(\varphi^-_N(\omega) \partial_x^{m+1}\varphi_N^-(\omega)) 
$$$$-\pi_{>N} (\partial_x^j \varphi^-_N
(\omega)\partial^{m-j+1}_x \varphi^-_N(\omega))\pi_{>N} 
(\varphi_N^+(\omega) \partial_x^{m+1}\varphi^+_N(\omega))].$$
Hence it is sufficient to show that
\begin{equation}\label{esca}
\limsup_{N\rightarrow \infty} \Big \|\int \pi_{>N} (\partial_x^j 
\varphi_N^+(\omega)\partial_x^{m-j+1} \varphi^+_N(\omega))
\pi_{>N} (\varphi_N^-(\omega) 
\partial_x^{m+1}\varphi_N^-(\omega))
\Big \|_{L^q_\omega}=0\end{equation}
$$\forall j=1,...,m.$$
Indeed the most delicate cases are 
$j=1, m$
(that in turn can be treated in a similar way). First we focus 
on \eqref{esca} for $j=m$. More precisely we have to prove
$$\limsup_{N\rightarrow \infty} \Big \|\int \pi_{>N} 
(\partial^m_x \varphi_N^+(\omega)\partial_x 
\varphi_N^+(\omega) )
\pi_{>N} (\varphi_N^-(\omega)\partial_x^{m+1}
\varphi_N^-(\omega)) dx \Big \|_{L^q_\omega}=0.$$
By replacing the random vector $\varphi(\omega)$ by its Fourier
randomized series we get:
$$\int \pi_{>N} (\partial^m_x \varphi^+_N(\omega)
\partial_x \varphi^+_N(\omega))
\pi_{>N} (\varphi_N^-(\omega) \partial_x^{m+1}
\varphi_N^-(\omega)) dx= $$
$$\sum_{\substack{0<|j_1|, |j_2|, |j_3|, |j_4|\leq N\\
j_1, j_2>0, j_3, j_4<0\\
|j_1+j_2|>N\\j_1+j_2+j_3+j_4=0}} \frac{\varphi_{j_1}(\omega)}{|j_1|} 
\frac{\varphi_{j_2}(\omega)}{|j_2|^m}
\frac{\varphi_{j_3}(\omega)}{|j_3|^{m+1}} \varphi_{j_4}(\omega).$$
Hence by the Minkowski inequality we get:
$$\limsup_{N\rightarrow \infty}  \Big \|\int \pi_{>N} 
(\partial^m_x \varphi_N^+(\omega)\partial_x 
\varphi_N^+(\omega) )
\pi_{>N} (\varphi_N^-(\omega)\partial_x^{m+1}
\varphi_N^-(\omega)) dx \Big \|_{L^q_\omega}
$$$$\leq C \limsup_{N\rightarrow \infty}
\Big (\sum_{0<|j_3|\leq N} \frac 1{|j_3|^{m+1}}
\Big) \Big (\sum_{\substack{0<|j_1|, |j_2|, |j_3|\leq N\\ |j_1+j_2|>N}}
\frac{1}{|j_1||j_2|^m}\Big )=O\Big (
\frac{\ln N}{N}\Big )$$
where we have used Lemma \ref{prod}.\\
Concerning the estimate \eqref{esca} for $j=1$ 
we can argue as above and we are reduced to prove that
$$\lim_{N\rightarrow \infty}
\sum_{\substack{0<|j_1|, |j_2|, |j_3|, |j_4|\leq N\\
j_1, j_2>0, j_3, j_4<0\\
|j_1+j_2|>N}} \frac 1{|j_1|^m 
|j_2||j_3|^{m+1}}=0$$
that follows by Lemma \ref{prod}.
The estimate for $III_N(\omega)$ is similar 
to the one of $II_N(\omega)$ provided that
\eqref{intpar22} is used instead of \eqref{intpar21}.
\\
\\
{\em Fourth case: $p(u)= u(\partial_x^m H u) (\partial_x^{m+1}u),
Hu\partial_x^m H u 
\partial_x^{m+1} Hu dx,
H u\partial_x^m u \partial_x^{m+1} u$}
\\
\\
They can be treated as in the third case provided that 
\eqref{mire1}, \eqref{mire2} (resp. 
\eqref{intpar4567},\eqref{intpar4567prime} 
and \eqref{intpar4567qua},\eqref{intpar4567quaprime}) are used instead of
\eqref{intpar21} and \eqref{intpar22}.

\hfill$\Box$

\begin{lem}\label{p3}
Let $m\geq 2$ be an integer and $p(u)\in {\mathcal P}_3(u)$ such that
$\tilde p(u)=\partial_x^\alpha u\partial_x^\beta 
u \partial_x^{\gamma}u$
with $$\alpha+\beta +\gamma=2m+1, 0\leq \alpha\leq 
\beta\leq \gamma \hbox{ and } \max
\{\alpha, \beta, \gamma\}\leq m.$$
Then we have 
$$\lim_{N\rightarrow \infty} \Big \|\int p^*_N(\pi_N 
u) dx \Big \|_{L^q(d\mu_{m+1})}=0, \hbox{ } \forall q\in [1,\infty).$$   
\end{lem}
{\bf Proof.} We treat for simplicity the case 
$p=\partial_x^\alpha u\partial_x^\beta u \partial_x^{\gamma}u$
(the general case can be studied with a similar argument).
Hence we get
$$p^*_N(\varphi_N(\omega))=I_N(\omega) + II_N(\omega)+III_N(\omega)$$
where
$$I_N(\omega)= \int \partial_x^\alpha(\pi_{>N} 
(\varphi_N(\omega) \partial_x \varphi_N(\omega)))
\partial_x^\beta \varphi_N(\omega) \partial_x^{\gamma} \varphi_N(\omega) dx;
$$
$$II_N(\omega)= \int 
\partial_x^\alpha \varphi_N(\omega) \partial_x^\beta 
(\pi_{>N} (\varphi_N(\omega)\partial_x \varphi_N(\omega))) 
\partial_x^{\gamma} \varphi_N(\omega) dx;
$$$$III_N(\omega)= \int \partial_x^\alpha 
\varphi_N(\omega)\partial_x^\beta  
\varphi_N(\omega)\partial^{\gamma}_x
(\pi_{>N} (\varphi_N(\omega)
\partial_x \varphi_N(\omega))) dx.$$ 
We shall prove that
$$\lim_{N\rightarrow \infty} \|I_N(\omega)\|_{L^q_\omega}=0$$
(and in a similar way we can treat $II_N(\omega)$ 
and $III_N(\omega)$). By 
the Leibnitz formula 
it is sufficient to prove 
$$
\lim_{N\rightarrow \infty}
\Big \|\int \pi_{>N} (\partial_x^j \varphi_N
(\omega) \partial_x^{\alpha-j+1}\varphi_N(\omega))
\partial_x^\beta \varphi_N(\omega) \partial_x^{\gamma} 
\varphi_N(\omega) dx\Big \|_{L^q_\omega}=0$$
$$\forall j=0,...,\alpha.$$
We shall treat the case $j=0$ and all the other 
cases can be treated in a similar way.
More precisely we shall prove that
$$
\lim_{N\rightarrow \infty}
\Big \|\int \pi_{>N} (\varphi_N(\omega) 
\partial_x^{\alpha+1}\varphi_N(\omega))
\partial_x^\beta \varphi_N(\omega) \partial_x^{\gamma} 
\varphi_N(\omega) dx\Big \|_{L^q_\omega}=0.$$
Notice that we have 
$$\int \pi_{>N} (\varphi_N(\omega) 
\partial_x^{\alpha+1}\varphi_N(\omega))
\partial_x^\beta \varphi_N(\omega) \partial_x^{\gamma} 
\varphi_N(\omega) dx$$$$=
\sum_{\substack{|j_1|, |j_2|, |j_3|, |j_4|\in (0,N],\\
|j_1+j_2|>N\\j_1+j_2+j_3+j_4=0}} \frac{\varphi_{j_1}(\omega)}
{|j_1|^{m+1}} \frac{\varphi_{j_2}(\omega)}{|j_2|^{m-\alpha}}
\frac{\varphi_{j_3}(\omega)}{|j_3|^{m+1-\beta}} 
\frac{\varphi_{j_4}(\omega)}{|j_4|^{m+1-\gamma}}$$
and hence
by using the triangular inequality we get
$$\Big \|\int \pi_{>N} (\varphi_N(\omega) 
\partial_x^{\alpha+1}\varphi_N(\omega))
\partial_x^\beta \varphi_N(\omega) \partial_x^{\gamma} 
\varphi_N(\omega) dx\Big \|_{L^q_\omega}$$$$\leq C
\sum_{\substack{|j_1|, |j_2|, |j_3|, |j_4|\in (0,N],\\
|j_1+j_2|>N\\j_1+j_2+j_3+j_4=0}} \frac 1{|j_1|^{m+1}
|j_2|^{m-\alpha} 
|j_3|^{m+1-\beta} |j_4|^{m+1-\gamma}}$$
Next we consider two possible cases:
\\
\\
{\em First subcase: $\alpha=1$, $\beta=\gamma=m$} 
\\
\\
In this case
we get
$$\Big \|\int (\pi_{>N} \varphi_N(\omega) 
\partial_x^{\alpha+1}\varphi_N(\omega))
\partial_x^\beta \varphi_N(\omega) 
\partial_x^{\gamma} \varphi_N(\omega) dx\Big \|_{L^q_\omega}$$
$$\leq C
\sum_{\substack{|j_1|, |j_2|, |j_4|\in (0,N],\\|j_1+j_2|>N}} 
\frac 1{|j_1|^{m+1}|j_2|^{m-1} |j_4|}$$
$$\leq \Big(\sum_{0<|j_4|\leq N} 
\frac 1{|j_4|}\Big)\Big(\sum_{\substack{0<|j_1|, |j_2|\leq N,\\ |j_1+j_2|>N}}
\frac 1{|j_1|^{m+1}|j_2|^{m-1}}\Big )
=O\Big (\frac{\ln ^2 N}{N} \Big) $$
where we have used Lemma \ref{prod}. 
\\
\\
{\em Second subcase: $\alpha\leq \beta=\gamma< m$} 
\\
\\
In this case we get
$$\Big \|\int (\pi_{>N} \varphi_N(\omega) 
\partial_x^{\alpha+1}\varphi_N(\omega))
\partial_x^\beta \varphi_N(\omega) \partial_x^{\gamma} 
\varphi_N(\omega) dx\Big \|_{L^q_\omega}$$$$\leq C
\sum_{\substack{|j_1|, |j_2|, |j_3|, |j_4|\in (0,N],\\
|j_1+j_2|>N}} \frac 1{|j_1|^{m+1}|j_2| 
|j_3|^{2}}$$
$$\leq C \Big(\sum_{0<|j_4|\leq N} 
\frac 1{|j_3|^2}\Big)\Big(\sum_{\substack{0<|j_1|, |j_2|\leq N,\\ |j_1+j_2|>N}}
\frac 1{|j_1|^{m+1}|j_2|}\Big )
=O\Big (\frac{\ln N}{N} \Big) $$
where we have used Lemma \ref{prod}.
\\
\\
{\em Third subcase: $\alpha\leq \beta<\gamma\leq m$} 
\\
\\
In this case we get
$$\Big \|\int (\pi_{>N} \varphi_N(\omega) 
\partial_x^{\alpha+1}\varphi_N(\omega))
\partial_x^\beta \varphi_N(\omega) \partial_x^{\gamma} 
\varphi_N(\omega) dx\Big \|_{L^q_\omega}$$
$$\leq C
\sum_{\substack{|j_1|, |j_2|, |j_3|, |j_4|\in (0,N],\\
|j_1+j_2|>N}} \frac 1{|j_1|^{m+1}|j_2| 
|j_3|^{2}}$$
and we can conclude as in the previous case. 

\hfill$\Box$

\begin{lem}\label{quartic}
Let $p(u)\in {\mathcal P}_4(u)$ be such that
$\tilde p(u)=\partial_x^{\alpha_1}u\partial_x^{\alpha_2}u 
\partial_x^{\alpha_3}u \partial_x^{\alpha_4}u$
with
$$\alpha_1\leq \alpha_2\leq \alpha_3\leq 
\alpha_4\leq m \hbox{ and } \alpha_1+
\alpha_2+\alpha_3+\alpha_4= 2m.$$
Then we have
$$\lim_{N\rightarrow \infty} \Big \|\int 
p^*_N(\pi_N u) dx\Big \|_{L^q(d\mu_{m+1})}=0, \hbox{ } \forall q\in [1,\infty).$$  
\end{lem}

{\bf Proof.} We shall treat the case 
$p=\partial_x^{\alpha_1}u\partial_x^{\alpha_2}u 
\partial_x^{\alpha_3}u \partial_x^{\alpha_4}u$. 
The general case follows in a similar way
(indeed our argument will be 
essentially based on the Minkowski inequality
and it is not affected in the case when $H$ appears in the expression 
of $p(u))$).
Arguing as in Lemma \ref{p3} it is sufficient to prove
that $$\lim_{N\rightarrow \infty} \|I_N\|_{L^q_\omega}=0, 
\lim_{N\rightarrow \infty} \|II_N\|_{L^q_\omega}=0,
\lim_{N\rightarrow \infty} \|III_N\|_{L^q_\omega}=0,
\lim_{N\rightarrow \infty} \|IV_N\|_{L^q_\omega}=0$$ where
$$I_N=\int \partial_x^{\alpha_1} \pi_{>N} 
(\varphi_N(\omega) \partial_x \varphi_N(\omega))
\partial_x^{\alpha_2} \varphi_N(\omega) 
\partial_x^{\alpha_3} \varphi_N(\omega) 
\partial_x^{\alpha_4} \varphi_N(\omega) dx$$
$$II_N=\int \partial_x^{\alpha_1} 
\varphi_N(\omega) \partial_x^{\alpha_2} 
\pi_{>N} (\varphi_N(\omega) \partial_x \varphi_N(\omega))
\partial_x^{\alpha_3} \varphi_N(\omega) \partial_x^{\alpha_4} 
\varphi_N(\omega) dx$$
$$III_N=\int \partial_x^{\alpha_1} \varphi_N(\omega) 
\partial_x^{\alpha_2} \varphi_N(\omega) \partial_x^{\alpha_3} 
\pi_{>N} (\varphi_N(\omega) \partial_x \varphi_N(\omega))
\partial_x^{\alpha_4} \varphi_N(\omega) dx$$
$$IV_N= \int \partial_x^{\alpha_1} \varphi_N(\omega) 
\partial_x^{\alpha_2} \varphi_N (\omega)
\partial_x^{\alpha_3} \varphi_N(\omega) \partial_x^{\alpha_4} 
\pi_{>N} (\varphi_N(\omega) \partial_x 
\varphi_N(\omega)) dx.
$$
We shall treat for simplicity 
only the term $IV_N$ (the other terms 
can be treated in a similar way). Hence we shall prove that
$\lim_{N\rightarrow \infty} \|IV_N\|_{L^q_\omega}=0$.
By the Leibnitz rule it follows by the following estimates: 
$$\lim_{N\rightarrow \infty} \Big \|\int \partial_x^{\alpha_1} 
\varphi_N(\omega) \partial_x^{\alpha_2} \varphi_N(\omega) 
\partial_x^{\alpha_3} \varphi_N (\omega) 
\pi_{>N} (\partial_x^j \varphi_N(\omega) \partial_x^{\alpha_4-j+1} 
\varphi_N(\omega)) dx\Big \|_{L^q_\omega}=0$$
$$\forall j=0,...,\alpha_4.$$
We shall prove the estimate above for $j=0$
(all the other cases can be treated in a simpler way).
Hence we have to show
$$\lim_{N\rightarrow \infty} \Big \|\int \partial_x^{\alpha_1} 
\varphi_N (\omega)\partial_x^{\alpha_2} \varphi_N (\omega)
\partial_x^{\alpha_3} \varphi_N (\omega) 
\pi_{>N} (\varphi_N(\omega) \partial_x^{\alpha_4+1} 
\varphi_N(\omega)) dx\Big \|_{L^q_\omega}=0.$$
Notice that
$$\int \partial_x^{\alpha_1} \varphi_N (\omega)
\partial_x^{\alpha_2} \varphi_N(\omega) 
\partial_x^{\alpha_3} \varphi_N  (\omega)
\pi_{>N} (\varphi_N(\omega) \partial_x^{\alpha_4+1} 
\varphi_N(\omega)) dx$$$$= \sum_{\substack{|j_1|, |j_2|, |j_3|, |j_4|, |j_5| \in 
(0,N],\\|j_4+j_5|>N\\j_1+j_2+j_3+j_4+j_5=0}} 
\frac{\varphi_{j_1}(\omega)}{|j_1|^{m+1-\alpha_1}}
\frac{\varphi_{j_2}(\omega)}{|j_2|^{m+1-\alpha_2}}\frac{\varphi_{j_3}(\omega)}
{|j_3|^{m+1-\alpha_3}}
\frac{\varphi_{j_4}(\omega)}{|j_4|^{m+1}}
\frac{\varphi_{j_5}(\omega)}{|j_5|^{m-\alpha_4}}$$
and hence by the Minkowski inequality
and Lemma \ref{orthspa} we get
$$\Big \|\int \partial_x^{\alpha_1} \varphi_N 
\partial_x^{\alpha_2} \varphi_N 
\partial_x^{\alpha_3} \varphi_N  
\pi_{>N} (\varphi_N \partial_x^{\alpha_4+1} 
\varphi_N) dx\Big \|_{L^q_\omega}$$$$\leq C
\sum_{\substack{|j_1|, |j_2|, |j_3|, |j_4|, |j_5| \in (0,N],\\|j_1+j_2+j_3|>N}}  
\frac{1}{|j_1|^2|j_2|^2|j_3||j_4|}=O\left (\frac{\ln^2 N}N\right )
$$
where we have used the fact that by assumption necessarily $\alpha_1, \alpha_2<m$.

\hfill$\Box$

\begin{lem}\label{multi}
Let $p_j(u)\in {\mathcal P}_j(u)$ with $j\geq 5$ be such that
$\tilde p_j(u)=\prod_{k=1}^j\partial_x^{\alpha_k}u$
where
$$\alpha_1\leq...\leq \alpha_j\leq m \hbox{ and }
\sum_{k=1}^j \alpha_k \leq 2m-1.$$
Then
$$\lim_{N\rightarrow \infty} \Big \|\int 
p^*_N(\pi_N u) dx\Big \|_{L^q(d\mu_{m+1})}=0, \hbox{ } \forall q\in [1,\infty).$$  
\end{lem}

{\bf Proof.} It follows as Lemma \ref{quartic}.

\hfill$\Box$

{\bf Proof of Theorem \ref{invariance}}
It follows by combining Lemma \ref{p3sing}, 
\ref{p3}, \ref{quartic}, \ref{multi}
with Proposition \ref{defG}.
\\

{\bf Acknowledgments.} We are grateful to the referee for the valuable remarks
and to Farinaz Wigmans for improving the presentation of the paper. The first author is partially supported by a grant from the European Research Council. 

\end{document}